\RequirePackage{fix-cm}
\documentclass[smallcondensed]{svjour3}     
\smartqed  
\usepackage{graphicx}
\usepackage{epic,eepic,bxeepic}
\usepackage{hhline}
\usepackage{amsmath}
\usepackage{amsfonts}
\usepackage{epstopdf}
\usepackage{multirow}
\usepackage{scalerel}
\usepackage{color}
\begin{document}

\title{Multigrid solvers for multipoint flux approximations of the Darcy problem on rough quadrilateral grids}

\titlerunning{Multigrid solvers for multipoint flux approximations of the Darcy problem}        

\author{Andr\'es Arrar\'as  \and Francisco J. Gaspar \and Laura Portero \and Carmen Rodrigo}


\institute{A. Arrar\'as \at
InaMat$^2$, Departamento de Estad\'{\i}stica, Inform\'atica y Matem\'aticas, Universidad P\'ublica de Navarra, Edificio de Las Encinas, Campus de Arrosad\'ia, 31006 Pamplona, Spain \\
 \email{andres.arraras@unavarra.es}           
\and
F.J. Gaspar \at
IUMA, Departamento de Matem\'atica Aplicada, Universidad de Zaragoza, Pedro Cerbuna 12, 50009 Zaragoza, Spain\\
\email{fjgaspar@unizar.es}
\and 
L. Portero \at 
InaMat$^2$, Departamento de Estad\'{\i}stica, Inform\'atica y Matem\'aticas, Universidad P\'ublica de Navarra, Edificio de Las Encinas, Campus de Arrosad\'ia, 31006 Pamplona, Spain \\
\email{laura.portero@unavarra.es}
\and 
C. Rodrigo \at
IUMA, Departamento de Matem\'atica Aplicada, Universidad de Zaragoza, Pedro Cerbuna 12, 50009 Zaragoza, Spain\\
\email{carmenr@unizar.es}
}

\date{Received: date / Accepted: date}

\maketitle

\begin{abstract}
In this work, an efficient blackbox-type multigrid method is proposed for solving multipoint flux approximations of the Darcy problem on logically rectangular grids. The approach is based on a cell-centered multigrid algorithm, which combines a piecewise constant interpolation and the restriction operator by Wesseling/Khalil with a line-wise relaxation procedure. A local Fourier analysis is performed for the case of a Cartesian uniform grid. The method shows a robust convergence for different full tensor coefficient problems and several rough quadrilateral grids.
\end{abstract}

\section{Introduction}
\label{sec:intro}

Single-phase flow in porous media is governed, under saturation conditions, by the first-order system
\begin{subequations}
	\label{cont:problem}
	\renewcommand{\theequation}{\theparentequation\alph{equation}}
	\begin{align}
	\mathbf{u}&=-K\,\nabla p,&&\hbox{in } \Omega,\label{cont:problem:a}\\
	\nabla\cdot \mathbf{u}&=f,&&\hbox{in } \Omega,\label{cont:problem:b}\\
	p&=g,&&\hbox{on } \Gamma_D,\label{cont:problem:c}\\
	\mathbf{u}\cdot\mathbf{n}&=0,&&\hbox{on } \Gamma_N.\label{cont:problem:d}
	\end{align}
\end{subequations}
Equations \eqref{cont:problem:a} and \eqref{cont:problem:b} represent Darcy's law and the continuity equation, respectively, and conform the so-called mixed formulation of the Darcy problem. Here, $\Omega$ is assumed to be a polygonal domain, whose boundary $\partial\Omega$ is divided into two non-overlapping subdomains, $\Gamma_D$ and $\Gamma_N$, in which Dirichlet and Neumann boundary conditions are imposed, respectively. Formally, $\partial\Omega=\overline{\Gamma}_D\cup\overline{\Gamma}_N$ such that $\Gamma_D\cap\Gamma_N=\emptyset$. The unknowns $p=p(\mathbf{x})$ and $\mathbf{u}=\mathbf{u}(\mathbf{x})$ stand for the fluid pressure and the Darcy velocity, respectively, and $K=K(\mathbf{x})$ is a $2\times 2$ symmetric positive definite tensor that represents the rock permeability divided by the fluid viscosity. This tensor is further assumed to satisfy
$$k_1\,\mathbf{\xi}^T\mathbf{\xi}\leq
\mathbf{\xi}^{T}K(\mathbf{x})\,\mathbf{\xi}\leq
k_2\,\mathbf{\xi}^T\mathbf{\xi},\qquad\quad\mathbf{x}\in\Omega,\,\mathbf{\xi}\neq\mathbf{0}\in\mathbb{R}^2,$$
where $0<k_1\leq k_2<\infty$. The problem is completed by providing data functions $f=f(\mathbf{x})$ and $g=g(\mathbf{x})$. Note that $\mathbf{n}$ denotes the outward unit normal on $\partial\Omega$. For the sake of simplicity, the gravitational term usually involved in Darcy's law has been neglected here. This simplification is acceptable since, in order to study the performance of the multigrid solver on the discretized problem, we are only concerned with the structure of the matrix and not with that of the right hand side.

Locally conservative discretization methods have been successfully applied to approximate the solution to problem \eqref{cont:problem}. In this framework, it is worth mentioning several techniques, such as mixed finite element (MFE) methods \cite{bre:for:91}, control volume MFE schemes \cite{cai:jon:mcc:rus:97,rus:whe:yot:07}, support-operator methods \cite{hym:sha:ste:97,sha:ste:96}, or the related mimetic finite difference methods \cite{bre:buf:lip:09,hym:mor:sha:ste:02,lip:man:sha:14}. These methods can all be formulated on irregular meshes and accurately handle anisotropic discontinuous tensors $K$. Moreover, they are designed to preserve the continuity of normal velocities across interelement edges. Equivalence relationships among them have been established in \cite{kla:rus:04}. The common drawback shared by all these methods is the need to solve an algebraic system with saddle point structure, which is known to be indefinite.

Several alternatives have been proposed to circumvent this problem, mainly based on a suitable choice of MFE spaces and the introduction of certain quadrature rules to approximate the vector-vector inner products. These strategies permit to construct an approximate Schur complement of the original system by locally eliminating the velocity unknowns, thus yielding a cell-centered finite difference scheme for the pressures. Examples of this idea on triangular meshes can be found in \cite{bar:mai:oud:96}, for the Laplace equation; in \cite{bre:for:mar:06}, for the Poisson equation; or in \cite{mic:sac:sal:01}, for a reaction--diffusion problem involving a diagonal tensor coefficient. Extensions to the Darcy system are proposed in \cite{rus:whe:83,wei:whe:88}, for diagonal tensors $K$ on rectangular meshes; or in \cite{arb:daw:kee:whe:yot:98}, for continuous full tensors $K$ on logically rectangular meshes obtained from smooth mappings of rectangular grids. To conclude, it should also be mentioned that the introduction of interelement Lagrange multipliers giving rise to the so-called mixed hybrid finite element methods \cite{arn:bre:85} is also a classical alternative for avoiding saddle point problems.\footnote{In this case, a multigrid solver can be applied by using the equivalence with the nonconforming Crouzeix--Raviart element (see, e.g., \cite{bau:kna:04,bra:ver:90}).}

In this work, we will focus on the so-called multipoint flux approximation (MPFA) methods, which combine the main advantages of the previous techniques. These methods were originally formulated as a natural extension of the well-known two-point flux approximation schemes. As such, they are defined to be control volume schemes which consider more than two pressure values to approximate the velocity at each cell edge. They have been formulated on unstructured triangular meshes \cite{aav:bar:boe:man:98a,aav:bar:boe:man:98b,edw:02}, and on logically rectangular meshes composed of quadrilateral elements \cite{aav:02,aav:bar:boe:man:96,edw:rog:98}. Following the terminology from \cite{aav:eig:kla:whe:yot:07}, these methods can be defined either in a reference space or in the physical space. Henceforth, we will describe the main properties of both approaches for quadrilateral elements. Remarkably, the triangular case has been analyzed in \cite{kla:rad:eig:08} using an equivalent MFE formulation, and a relevant extension for solving Richards' equation has been further proposed. 

MPFA schemes in a reference space can be formulated as MFE methods using the trapezoidal quadrature rule, with a particular choice of finite element spaces and evaluation points: the method proposed in \cite{kla:win:06a} is based on the so-called broken Raviart--Thomas element, $\mathcal{RT}_{1/2}$, with an evaluation of the tensor at the midpoint of the reference element; in turn, that proposed in \cite{whe:yot:06} is based on the lowest order Brezzi--Douglas--Marini element, $\mathcal{BDM}_1$, with an evaluation of the tensor at the corners of the reference element.   Note that, in both cases, the degrees of freedom of the corresponding MFE method show a one-to-one correspondence with the unknowns of the MPFA scheme. In addition, they both result in a cell-centered pressure system with a symmetric positive definite matrix. Numerical evidence reveals that the second variant is more accurate than the first one when computing normal velocities on grids composed of $\mathcal{O}(h^2)$-perturbations of parallelograms \cite{aav:eig:kla:whe:yot:07}. In turn, a loss of convergence is observed for both methods on $\mathcal{O}(h)$-perturbed meshes. In this case, the physical space MPFA method permits to recover a first-order convergence rate for the normal velocities. The price to pay is a loss of symmetry of the system matrix unless the mesh is composed of parallelograms.

MPFA methods in the physical space can also be formulated as MFE schemes, based on the $\mathcal{RT}_{1/2}$ or the $\mathcal{BDM}_1$ element, using a non-symmetric variant of the trapezoidal quadrature rule. In the former case, the quadrature rule involves the mean value of the inverse tensor $K^{-1}$ on each element $E$, denoted by $\overline{K_E^{-1}}$ \cite{kla:win:06}; in the latter, it considers the inverse of the mean value of $K$ on $E$, that is, $\overline{K}_E^{-1}$ \cite{whe:xue:yot:12a}. The solvability of the resulting discrete problem is only guaranteed under technical restrictions on the element geometry and the anisotropy of $K$ (see, e.g., formulas (3.31) and (3.32) in \cite{whe:xue:yot:12a}; cf. also \cite{lip:sha:yot:09}). Remarkably, when the mesh is composed of parallelograms and $K$ is piecewise constant on each element, the physical space MPFA methods are equivalent to their corresponding reference space counterparts.

In the present paper, we deal with the efficient solution of the large systems of equations arising from the reference and physical space MPFA methods based on the $\mathcal{BDM}_1$ finite element spaces. According to \cite{whe:xue:yot:12a,whe:yot:06}, we will henceforth rename such methods the symmetric and non-symmetric multipoint flux mixed finite element (MFMFE) schemes, respectively. Broadly speaking, the design of solvers for MPFA methods has been seldom considered in the literature. This is somewhat surprising due to the wide use of these schemes in the numerical simulation of porous media problems. A two-level domain decomposition algorithm for an MPFA finite volume discretization of three-dimensional flow in anisotropic heterogeneous porous media is presented in \cite{Iryna_two_level}. In this interesting work, additive and multiplicative Schwarz iterative methods are implemented as smoothers, and the coarse scale operator is obtained from numerical upscaling. To the best of our knowledge, however, a robust multigrid solver for multipoint flux discretizations has never been proposed so far. This work is intended to close this gap, presenting the first multigrid method for solving multipoint flux approximations, which can easily handle difficult diffusion problems on rough grids and with coefficient jumps not aligned with the mesh.

It is well known that present multigrid methods are among the most advanced techniques for solving large linear systems arising from the discretization of partial differential equations (PDEs). There are mainly two different approaches to multigrid: geometric multigrid methods (GMG), for which a hierarchy of grids has to be defined and the inter-grid transfer operations are based on geometric principles, and algebraic multigrid methods (AMG), which construct the coarse levels automatically from the system matrix on the target grid in an algebraic way. In the context of AMG methods, two prevailing schemes have proved their use for multiple engineering problems, namely: algebraic multigrid and aggregation-based multigrid methods \cite{Braess1995,Brezina,stuben2001introduction,Van2001,Van1996}. The origin of algebraic methods may be found in the early days of multigrid, when blackbox multigrid (BoxMG) with operator-dependent transfer operators and Galerkin coarse grid approximation was proposed for structured vertex-centered Cartesian grids \cite{Alcouffe:1981:MGM,dendy0,DENDY1983261}. This approach can be considered as a predecessor of classical AMG. The aggregation-based multigrid methods, whose origin can be found in \cite{Van1996} (smoothed aggregation), may be related to the cell-centered multigrid methods proposed in \cite{MG2,WESSELING198885}. In \cite{WESSELING198885}, it was shown that constant (i.e., operator-independent) transfer operators, in combination with Galerkin coarse grid discretization, provided highly efficient multigrid results for cell-centered discretizations of elliptic PDEs including jumping coefficients. In addition, more robustness can be achieved in some cases by using multigrid as a preconditioner of a Krylov subspace method; see, e.g., \cite{ash:fal:96}, where a parallel multigrid preconditioned conjugate gradient method is proposed for solving groundwater flow problems. Like AMG and BoxMG, the multigrid solver proposed here is based on the blackbox methodology, in which the user only provides the fine-grid discretization, the right-hand side and an initial guess for the solution, and the code automatically generates the hierarchy of operators on coarser levels.  

Logically rectangular grids are considered to take advantage of the recent trends in computer architectures: many-core and accelerated architectures that achieve their best performance when structured data can be used. On this type of architectures, the drawbacks of using unstructured data access and indirect addressing may be so significant that it becomes more efficient to consider larger, logically structured grids with regular data access. In particular, implementations of logically structured multigrid algorithms, such as BoxMG, have been shown to be 10 times faster than AMG for three-dimensional heterogeneous diffusion problems on structured grids \cite{Moulton2012}.

In the present work, a local Fourier analysis (LFA) \cite{Bra77,Bra94} is performed for the case of Cartesian rectangular grids. This analysis is known to predict very accurately the convergence rates of GMG methods by considering the local character of the involved operators and neglecting the effect of boundary conditions. Specifically, the Toeplitz or multilevel Toeplitz structure of the matrix on an infinite grid permits its diagonalization by the matrix of Fourier modes. A detailed introduction to this analysis can be found in \cite{TOS01,Wie01}. The robustness of the proposed multigrid solver is shown for different test problems with full tensor coefficients on various rough logically rectangular grids.

The remaining of the paper is structured as follows. In Section \ref{sec:1}, we describe the considered MFMFE schemes for the Darcy system. Subsection \ref{sec:1.1} is devoted to deriving the stencil coefficients for homogeneous media and Cartesian uniform grids, which will be subsequently used in the LFA. In Section \ref{sec:2}, the cell-centered blackbox-type multigrid method is explained. This section also includes a description of the LFA, together with some predicting results of the convergence rates of the multigrid method. Section \ref{sec:3} shows the robustness of the proposed multigrid solver for problems with different permeability tensors on a variety of logically rectangular grids. Finally, Section \ref{sec:conclusions} provides some concluding remarks together with certain ideas for future research.

\section{Multipoint flux approximations of the Darcy system}
\label{sec:1}

The mixed variational formulation of the first-order system (\ref{cont:problem}) reads: \emph{Find} $(\mathbf{u},p)\in V\times W$ \emph{such that}
\begin{subequations}\label{weak:mixed:formulation}
	\renewcommand{\theequation}{\theparentequation\alph{equation}}
	\begin{align}
	(K^{-1}\mathbf{u},\mathbf{v})&=(p,\nabla\cdot\mathbf{v})-\langle g,\mathbf{v}\cdot\mathbf{n}\rangle_{\Gamma_D},&&\mathbf{v}\in V,\label{weak:mixed:formulation:a}\\[0.5ex]
	(\nabla\cdot\mathbf{u},w)&=(f,w),&&w\in W,\label{weak:mixed:formulation:b}
	\end{align}
\end{subequations}
where $(\cdot,\cdot)$ denotes the inner product in either $L^2(\Omega)$ or $(L^2(\Omega))^2$, and $\langle\cdot,\cdot\rangle_{\Gamma_D}$ represents the $L^2(\Gamma_D)$-inner product or duality pairing. Moreover, $W=L^2({\Omega})$ and $V=\left\{\mathbf{v}\in {H}(\mbox{div};\Omega):\mathbf{v}\cdot\mathbf{n}=0\ \mathrm{on}\ \Gamma_N\right\}$, with
$$
H(\mbox{div};G)=\{\mathbf{v}\in (L^2(G))^2:\nabla\cdot\mathbf{v}\in L^2(G)\},
$$
for any $G\subset\mathbb{R}^2$. It is well known that (\ref{weak:mixed:formulation}) has a unique solution \cite{bre:for:91}.
For the sake of simplicity, we will henceforth consider homogeneous Dirichlet boundary conditions along $\Gamma_D$, i.e., $g\equiv 0$.

Next, we introduce an MFMFE discretization for (\ref{weak:mixed:formulation}) based on the $\mathcal{BDM}_1$ spaces on quadrilateral elements \cite{bre:dou:mar:85,bre:for:91}. As we will see below, the method further considers suitable quadrature rules which permit to eliminate the velocity unknowns, thus yielding a cell-centered  scheme for the pressure.

Let $\mathcal{T}_h$ be a logically rectangular partition of $\Omega$ into convex quadrilaterals, where $h=\max_{E\in\mathcal{T}_h}\mathrm{diam}(E)$. We denote by $\hat{E}$ the unit square with vertices $\hat{\mathbf{r}}_1=(0,0)^T$, $\hat{\mathbf{r}}_2=(1,0)^T$, $\hat{\mathbf{r}}_3=(1,1)^T$ and $\hat{\mathbf{r}}_4=(0,1)^T$, and introduce a family of bijective bilinear mappings $\{F_E\}_{E\in\mathcal{T}_h}$ such that $F_E(\hat{E})=E$. 
In particular, given a physical element $E\in\mathcal{T}_h$ with vertices $\mathbf{r}_i=(x_i,y_i)^T$, for $i=1,2,3,4$, $F_E$ is defined as
$$F_E(\hat{\mathbf{r}})=\mathbf{r}_1+\mathbf{r}_{21}\hat{x}+\mathbf{r}_{41}\hat{y}+(\mathbf{r}_{34}-\mathbf{r}_{21})\hat{x}\hat{y},$$
where $\mathbf{r}_{ij}=\mathbf{r}_i-\mathbf{r}_j$.
We further define, for each mapping $F_E$, the Jacobian matrix $DF_E$ and its determinant $J_E=|\det(DF_E)|$. The outward unit vectors normal to the edges $\hat{e}\subset\partial\hat{E}$ and $e\subset\partial E$ are denoted by $\hat{\mathbf{n}}_{\hat{e}}$ and $\mathbf{n}_e$, respectively.

For later use, a given partition $\mathcal{T}_h$ will be called an $\mathcal{O}(h^2)$-perturbed mesh if it is composed of $h^2$-parallelograms, that is, quadrilaterals whose vertices satisfy the condition $|\mathbf{r}_{34}-\mathbf{r}_{21}|_{\mathbb{R}^2}\leq C h^2$. This kind of elements can be obtained asymptotically by uniform refinements of a general quadrilateral grid.

\begin{figure}[t]
	\begin{center}
		\vspace*{0.45cm}\hspace*{0.15cm}\unitlength=0.035cm\begin{picture}(310,100)
		\allinethickness{0.75pt}
		\put(30,30){\line(1,0){60}}\put(90,30){\line(0,1){60}}\put(90,90){\line(-1,0){60}}\put(30,90){\line(0,-1){60}}
		\put(57,57){$\times$}
		\put(56,36){\footnotesize{${\hat{E}}$}}
		\put(35,30){\circle*{4}}
		\put(84.7,30){\circle*{4}}
		\put(89.9,36){\circle*{4}}
		\put(89.9,84){\circle*{4}}
		\put(84.7,90){\circle*{4}}
		\put(35.7,90){\circle*{4}}
		\put(30,84){\circle*{4}}
		\put(30,36){\circle*{4}}
		\put(22,23.5){\footnotesize{$\hat{\mathbf{r}}_1$}}
		\put(92,23.5){\footnotesize{$\hat{\mathbf{r}}_2$}}
		\put(92,91.5){\footnotesize{$\hat{\mathbf{r}}_3$}}
		\put(18,91.5){\footnotesize{$\hat{\mathbf{r}}_4$}}
		\put(122,53){\vector(1,0){60}}
		\put(146,57){\footnotesize{$F_{E}$}} \allinethickness{0.75pt}
		\path(210,40)(270,20)(284,66)(256,86)(210,40)
		\put(250,49){$\times$}
		\put(240,36.5){\footnotesize{${E}$}}
		\put(214.7,45){\circle*{4}}
		\put(217.7,37.6){\circle*{4}}
		\put(263.7,22.1){\circle*{4}}
		\put(271.7,26.1){\circle*{4}}
		\put(282,60){\circle*{4}}
		\put(279.7,69.2){\circle*{4}}
		\put(261.2,82.3){\circle*{4}}
		\put(251.3,81.5){\circle*{4}}
		\put(201,36.5){\footnotesize{$\mathbf{r}_1$}}
		\put(270,12.5){\footnotesize{$\mathbf{r}_2$}}
		\put(287,64.5){\footnotesize{$\mathbf{r}_3$}}
		\put(252,89.5){\footnotesize{$\mathbf{r}_4$}}
		\end{picture}
		\caption{Degrees of freedom for the $\mathcal{BDM}_1$ spaces on quadrilaterals; crosses and circles denote pressure and velocity degrees of freedom, respectively. For the sake of clarity, the velocity degrees of freedom are drawn at a distance from the corresponding vertex.}
		\label{fig:velocity:dof}
	\end{center}
\end{figure}
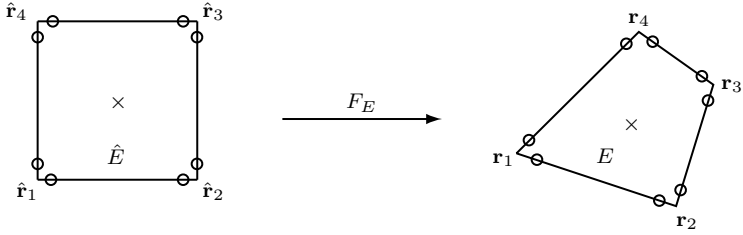

Let $\hat{V}(\hat{E})$ and $\hat{W}(\hat{E})$ be the $\mathcal{BDM}_1$ finite element spaces on the reference element $\hat{E}$, i.e.,
\begin{equation*}
\hat{V}(\hat{E})=\begin{bmatrix}
\alpha_1+\beta_1\hat{x}+\gamma_1\hat{y}+r\hat{x}^2+2s\hat{x}\hat{y}\\
\alpha_2+\beta_2\hat{x}+\gamma_2\hat{y}-2r\hat{x}\hat{y}-s\hat{y}^2\\
\end{bmatrix},
\qquad
\hat{W}(\hat{E})=\mathbb{P}_0(\hat{E}),\nonumber
\end{equation*}
where $\alpha_1$, $\alpha_2$, $\beta_1$, $\beta_2$, $\gamma_1$, $\gamma_2$, $r$ and $s$ are real constants. Note that $\hat{\nabla}\cdot\hat{V}(\hat{E})=\hat{W}(\hat{E})$ and, on any edge $\hat{e}\subset\partial\hat{E}$, $\hat{\mathbf{v}}\in\hat{V}(\hat{E})$ is such that $\hat{\mathbf{v}}\cdot\hat{\mathbf{n}}_{\hat{e}}\in\mathbb{P}_1(\hat{e})$. The degrees of freedom for $\hat{\mathbf{v}}$ are chosen to be the values of $\hat{\mathbf{v}}\cdot\hat{\mathbf{n}}_{\hat{e}}$ at the vertices of each edge $\hat{e}$ (see Figure \ref{fig:velocity:dof}). 

Any approximation to $H(\mbox{div};E)$ must preserve the continuity of the normal components of the velocity vectors across the element edges. In order to satisfy this requirement, the space $V(E)$ associated to an element  $E\in\mathcal{T}_h$ will be defined by using the Piola transformation \cite{bre:for:91}, i.e.,
\begin{equation*}
\mathbf{v}\leftrightarrow\hat{\mathbf{v}}:\mathbf{v}=(J_E^{-1}{DF}_E\,\hat{\mathbf{v}})\circ{F}_E^{-1}.
\end{equation*}
On the other hand, the transformation 
\begin{equation*}w\leftrightarrow\hat{w}: w = \hat{w}\circ{F}_E^{-1}
\end{equation*}
will be used for the definition of the local space $W(E)$. Finally, the global MFE spaces $V_h\times W_h\subset V\times W$ on $\mathcal{T}_h$ are given by
\begin{align}
V_h&=\{\mathbf{v}\in V: \mathbf{v}|_E\leftrightarrow\hat{\mathbf{v}},\,\hat{\mathbf{v}}\in \hat{V}(\hat{E})\,\,\forall\,E\in\mathcal{T}_h\},\nonumber\\
W_h&=\{w\in W: w|_E\leftrightarrow\hat{w},\,\hat{w}\in\hat{W}(\hat{E})\,\,\forall\,E\in\mathcal{T}_h\}.\nonumber
\end{align}

In the approximation of (\ref{weak:mixed:formulation:a}) by the MFE method described above, it is necessary to compute integrals of the form $(K^{-1}\,\mathbf{q},\mathbf{v})$, for $\mathbf{q},\mathbf{v}\in V_h$. In this setting, the MFMFE method is derived by considering suitable quadrature rules that allow for local velocity elimination. In particular, for any $\mathbf{q},\mathbf{v}\in V_h$, the global quadrature rule is defined element-wise as 
\begin{equation}\label{global:quadrature}
(K^{-1}\mathbf{q},\mathbf{v})_{Q}=\sum_{E\in\mathcal{T}_h}(K^{-1}\mathbf{q},\mathbf{v})_{Q,E}.
\end{equation}
The definition of the local quadrature rule depends on the type of spatial meshes under consideration. For $\mathcal{O}(h^2)$-perturbed meshes, the numerical integration on each element is defined by mapping to the reference element \cite{whe:yot:06}, i.e.,
\begin{equation}\label{quadr:symm}
(K^{-1}\,\mathbf{q},\mathbf{v})_{Q,E}=(\mathcal{K}_E^{-1}\hat{\mathbf{q}},\hat{\mathbf{v}})_{\hat{Q},\hat{E}}=
\frac{1}{4}\sum_{i=1}^{4}\mathcal{K}_E^{-1}(\hat{\mathbf{r}}_i)\hat{\mathbf{q}}(\hat{\mathbf{r}}_i)\cdot\hat{\mathbf{v}}(\hat{\mathbf{r}}_i),
\end{equation}
where
$$
\mathcal{K}_E^{-1}(\hat{\mathbf{x}})=\dfrac{1}{J_E(\hat{\mathbf{x}})}\,DF_E^T(\hat{\mathbf{x}})K^{-1}(F_E(\hat{\mathbf{x}}))\,DF_E(\hat{\mathbf{x}}).
$$
As for the case of highly distorted rough grids, the quadrature rule on each element is defined accordingly \cite{whe:xue:yot:12a}, that is,
\begin{equation}\label{quadr:non:symm}
(K^{-1}\mathbf{q},\mathbf{v})_{Q,E}=(\widetilde{\mathcal{K}}_E^{-1}\hat{\mathbf{q}},\hat{\mathbf{v}})_{\hat{Q},\hat{E}}=
\frac{1}{4}\sum_{i=1}^{4}\widetilde{\mathcal{K}}_E^{-1}(\hat{\mathbf{r}}_i)\hat{\mathbf{q}}(\hat{\mathbf{r}}_i)\cdot\hat{\mathbf{v}}(\hat{\mathbf{r}}_i),
\end{equation} 
where
$$
\widetilde{\mathcal{K}}_E^{-1}(\hat{\mathbf{x}})=\dfrac{1}{J_E(\hat{\mathbf{x}})}\,DF_E^T(\hat{\mathbf{x}}_c)\overline{K}_E^{-1}\,DF_E(\hat{\mathbf{x}}),
$$
with $\overline{K}_E$ being a constant matrix such that $\overline{K}_{ij,E}$ is the mean value of $K_{ij}$ on $E$, where $\overline{K}_{ij,E}$ and $K_{ij}$ are the elements on the $i$th row and $j$th column of matrices $\overline{K}_E$ and $K$, respectively. Furthermore, $\hat{\mathbf{x}}_{c}$ denotes the center of mass of $\hat{E}$. Note that, if $K$ is an element by element piecewise constant tensor and the spatial mesh consists of parallelograms, the expression (\ref{quadr:non:symm}) reduces to (\ref{quadr:symm}).

The MFMFE approximation to (\ref{weak:mixed:formulation}), considering homogeneous Dirichlet boundary conditions, is given by: \emph{Find} $(\mathbf{u}_h,p_h)\in {V_h}\times {W_h}$ \emph{such that}
\begin{subequations}\label{mfmfe:method}
	\renewcommand{\theequation}{\theparentequation\alph{equation}}
	\begin{align}
	(K^{-1}\mathbf{u}_h,\mathbf{v})_Q&=(p_h,\nabla\cdot\mathbf{v})
	,&&\mathbf{v}\in V_h,\label{mfmfe:method:a}\\[0.5ex]
	(\nabla\cdot\mathbf{u}_h,w)&=(f,w),&&w\in W_h.\label{mfmfe:method:b}
	\end{align}
\end{subequations}
If $(K^{-1}\cdot,\cdot)_Q$ in \eqref{mfmfe:method:a} is given by \eqref{global:quadrature} and \eqref{quadr:symm}, the MFMFE method is called symmetric; if it is given by \eqref{global:quadrature} and \eqref{quadr:non:symm}, it is referred to as non-symmetric. The well-posedness and convergence properties of both methods stemming from (\ref{mfmfe:method}) are studied in \cite{whe:yot:06} and \cite{whe:xue:yot:12a}, respectively. In the symmetric case, if the mesh is composed of $h^2$-parallelograms, the velocity and pressure variables are first-order convergent in the $L^2$-norm, the latter being second-order superconvergent at the cell centers in a discrete $L^2$-norm. In the non-symmetric case, the velocity variable is shown to be first-order convergent, either when compared to the projection of the true solution onto the space $V_h$ in the $L^2$-norm, or when considering the normal component in an edge-based norm; in turn, the pressure preserves the first-order optimal convergence in the $L^2$-norm. Note that the non-symmetric MFMFE scheme need not be applied on an $\mathcal{O}(h^2)$-perturbed mesh. The convergence properties of both methods will be illustrated by numerical experiments in the last section of the paper.

In the sequel, we describe how the MFMFE formulation (\ref{mfmfe:method}) can be reduced to a cell-centered finite difference scheme in the pressure variable. More precisely, we use (\ref{mfmfe:method:a}) to express the velocity degrees of freedom in terms of the pressure unknowns, and substitute them back into (\ref{mfmfe:method:b}) in order to obtain a linear system for the pressures.

In this framework, let us introduce the vector space $\mathcal{H}_v$ of discrete velocity functions $U_h\in\mathbb{R}^{2N_{\ell}}$, where $N_{\ell}$ is the number of edges in $\mathcal{T}_h$. The degrees of freedom for this space are located at the vertices of each edge. If $N_v$ denotes the number of vertices in $\mathcal{T}_h$, any $U_h\in\mathcal{H}_v$ is given by $U_h=(U_{h,1},U_{h,2},\ldots,U_{h,N_v})^T$, where $U_{h,i}\in\mathbb{R}^{\ell_i}$ and $\ell_{i}$ is the number of edges that share the $i$th vertex point, for $i=1,2,\ldots,N_v$. The component of $U_{h,i}$ associated to the $j$th edge $e_j$ is given by the volumetric flux $(\mathbf{u}_h\cdot\mathbf{n}_{e_j})(\mathbf{r}_i)|e_j|$, where $|e_j|$ denotes the length of $e_j$, for $j=1,2,\ldots,\ell_i$. On the other hand, we consider the vector space $\mathcal{H}_p$ of discrete pressure functions $P_{h}\in\mathbb{R}^{N_e}$, where $N_e$ stands for the number of elements in $\mathcal{T}_h$. In this case, the degrees of freedom are located at the cell centers. Any $P_h\in\mathcal{H}_p$ has the form $P_h=(P_{h,1},P_{h,2},\ldots,P_{h,N_e})^T$, where $P_{h,i}=p_h(\mathbf{x}_{c,i})$ and $\mathbf{x}_{c,i}$ is the coordinate vector of the $i$th cell center, for $i=1,2,\ldots,N_e$.

Let $\{\mathbf{v}_i\}_{i=1}^{L}$ and $\{w_j\}_{j=1}^{N_e}$ be the bases of the discrete spaces $V_h$ and $W_h$, respectively. The linear system stemming from (\ref{mfmfe:method:a})-(\ref{mfmfe:method:b}) is 
\begin{equation}\label{global:linear:system}
\begin{bmatrix}
A & B \\
B^T & 0
\end{bmatrix}
\begin{bmatrix}
{U}_h \\
P_h \\
\end{bmatrix}
=
\begin{bmatrix}
0\\
F_h \\
\end{bmatrix},
\end{equation}
where the matrices $A\in\mathbb{R}^{2N_{\ell}\times 2N_{\ell}}$ and $B\in\mathbb{R}^{2N_{\ell}\times N_e}$ are given by $(A)_{ij}=(K^{-1}\mathbf{v}_j,\mathbf{v}_i)_Q$ and
$(B)_{ij}=-(w_j,\nabla\cdot\mathbf{v}_i)$, respectively.
Finally, $F_h\in\mathbb{R}^{N_{e}}$ contains the contributions of the right-hand side $f$ and is defined element-wise as
$$(F_{h})_{E}=-\int_{E}f(\mathbf{x})\,d\mathbf{x}.$$ 

The use of the trapezoidal quadrature rule $(K^{-1}\cdot,\cdot)_Q$ permits to decouple the velocity degrees of freedom associated to a vertex from the rest of them. In consequence, the velocity mass matrix $A=\mathrm{diag}(A_{1},A_{2},\ldots,A_{N_v})$ has a block-diagonal structure, and each block $A_i\in\mathbb{R}^{\ell_{i}\times\ell_{i}}$ is a local matrix associated to the $i$th vertex point, for $i=1,2,\ldots,N_v$. Since $A$ is induced by the discrete bilinear form $(K^{-1}\cdot,\cdot)_Q$, the symmetry of the corresponding blocks $A_i$ depends on the local quadrature rule under consideration, namely: if $(K^{-1}\cdot,\cdot)_{Q,E}$ is given by (\ref{quadr:symm}), $A_i$ will be symmetric; otherwise, if it is given by (\ref{quadr:non:symm}), $A_i$ will be non-symmetric unless the mesh is composed of parallelograms and $K$ is constant on each element.

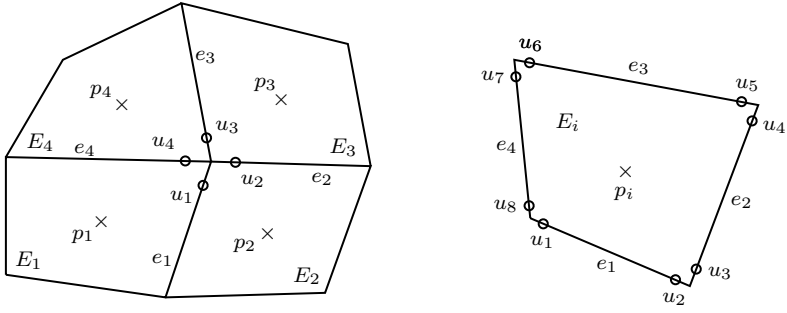
\begin{figure}[t]
	\begin{center}
		\label{fig:int:node:2}\unitlength=0.03cm
		\begin{picture}(400,170)
		\allinethickness{0.75pt}
		\path(30,30)(100,20)(120,80)(30,82)(30,30)\path(100,20)(170,22)(190,78)(120,80)\path(190,78)(180,132)(107,150)(120,80)\path(107,150)(55,125)(30,82)
		\put(116.5,69.5){\circle*{4}}
		\put(130.7,79.7){\circle*{4}}
		\put(118,90.5){\circle*{4}}
		\put(108.7,80.3){\circle*{4}}
		\put(102,63){{\footnotesize$u_{1}$}}
		\put(133,70){{\footnotesize$u_{2}$}}
		\put(122,94){{\footnotesize$u_{3}$}}
		\put(94,87){{\footnotesize$u_{4}$}}
		\put(94,35){\footnotesize{$e_{1}$}}
		\put(164,71){\footnotesize{$e_{2}$}}
		\put(113,125){\footnotesize{$e_{3}$}}
		\put(60,84){\footnotesize{$e_{4}$}}
		\put(68,51){$\times$}
		\put(141,46){$\times$}
		\put(147,105){$\times$}
		\put(77,103){$\times$}
		\put(59,46){\footnotesize{$p_{1}$}}
		\put(130,42){\footnotesize{$p_{2}$}}
		\put(139,113){\footnotesize{$p_{3}$}}
		\put(67,110){\footnotesize{$p_{4}$}}
		\put(34,34){\footnotesize{$E_{1}$}}
		\put(156,27){\footnotesize{$E_{2}$}}
		\put(172,84){\footnotesize{$E_{3}$}}
		\put(39,86){\footnotesize{$E_{4}$}}		
		\path(260,55)(330,25)(360,105)(253,125)(260,55)
		\put(298,73){$\times$}
		\put(297,65){\footnotesize{$p_{i}$}}
		\put(265.7,52.5){\circle*{4}}
		\put(323.7,27.8){\circle*{4}}
		\put(332.8,32.5){\circle*{4}}
		\put(357.3,98){\circle*{4}}
		\put(352.7,106.5){\circle*{4}}
		\put(259.7,123.7){\circle*{4}}
		\put(253.7,117.5){\circle*{4}}
		\put(259.5,60.5){\circle*{4}}
		\put(260,43){{\footnotesize$u_{1}$}}
		\put(318,18){{\footnotesize$u_{2}$}}
		\put(338,29.5){{\footnotesize$u_{3}$}}
		\put(362,94){{\footnotesize$u_{4}$}}
		\put(350,112.5){{\footnotesize$u_{5}$}}
		\put(255,130){{\footnotesize$u_{6}$}}
		\put(255,130){{\footnotesize$u_{6}$}}
		\put(238.5,116){{\footnotesize$u_{7}$}}
		\put(244,59){{\footnotesize$u_{8}$}}
		\put(271,95){\footnotesize{$E_{i}$}}
		\put(289,32.5){\footnotesize{$e_{1}$}}
		\put(348,61){\footnotesize{$e_{2}$}}
		\put(303,120){\footnotesize{$e_{3}$}}
		\put(245,85){\footnotesize{$e_{4}$}}
		\end{picture}
		\caption{Interaction of degrees of freedom in the MFMFE method.}
	\end{center}
\end{figure}

In particular, let us consider the $i$th interior vertex represented in Figure \ref{fig:int:node:2} (left). We denote by $u_j$, for $j=1,2,3,4$, the velocity degrees of freedom associated to this vertex. The corresponding basis functions of the space $V_h$ are denoted by $\mathbf{v}_j$, for $j=1,2,3,4$. Inserting $\mathbf{v}=\mathbf{v}_j$ into (\ref{mfmfe:method:a}), we obtain the local system corresponding to the $i$th interior vertex, i.e.,
\begin{equation}
\label{local:system:vertex}
A_i\begin{bmatrix}
u_1\,|e_1|\\u_2\,|e_2|\\u_3\,|e_3|\\u_4\,|e_4|\end{bmatrix}+\frac{1}{2}\begin{bmatrix}
p_2-p_1\\p_3-p_2\\p_3-p_4\\p_4-p_1\end{bmatrix}=\begin{bmatrix}
0\\0\\0\\0\end{bmatrix},
\end{equation}
where
\begin{equation}
\label{block:vertex:matrix}
A_i=\frac{1}{4}\begin{bmatrix}
\mathcal{N}_{11,E_1}^{-1}+\mathcal{N}_{11,E_2}^{-1} & \mathcal{N}_{12,E_2}^{-1} & 0 & \mathcal{N}_{12,E_1}^{-1}\\[1ex]
\mathcal{N}_{21,E_2}^{-1} & \mathcal{N}_{22,E_2}^{-1}+\mathcal{N}_{22,E_3}^{-1} &  \mathcal{N}_{21,E_3}^{-1} & 0 \\[1ex]
0 & \mathcal{N}_{12,E_3}^{-1} & \mathcal{N}_{11,E_3}^{-1}+\mathcal{N}_{11,E_4}^{-1} &  \mathcal{N}_{12,E_4}^{-1}\\[1ex]
\mathcal{N}_{21,E_1}^{-1} & 0 &  \mathcal{N}_{21,E_4}^{-1}  & \mathcal{N}_{22,E_4}^{-1}+\mathcal{N}_{22,E_1}^{-1} 
\end{bmatrix}.
\end{equation}
In this expression, using a generic element notation, $\mathcal{N}_{ij,E}^{-1}=\mathcal{K}_{ij,E}^{-1}$, if we consider the symmetric quadrature rule (\ref{quadr:symm}), and $\mathcal{N}_{ij,E}^{-1}=\widetilde{\mathcal{K}}_{ij,E}^{-1}$, for its non-symmetric counterpart (\ref{quadr:non:symm}) (see, e.g., \cite{whe:xue:yot:12a,whe:yot:06}).

On the other hand, let us consider the $i$th element represented in Figure \ref{fig:int:node:2} (right). We denote by $w_i$ the basis function of the space $W_h$ associated to this element. If we insert $w=w_i$ into (\ref{mfmfe:method:b}), we obtain the expression
\begin{equation}\label{discrete:div}
\frac{1}{2}\left((u_5+u_6)|e_3|-(u_1+u_2)|e_1|+(u_3+u_4)|e_2|-(u_7+u_8)|e_4|\right)=\int_{E_i}f(\mathbf{x}),
\end{equation}
where we consider the local notations introduced in Figure \ref{fig:int:node:2} (right). Therefore, the $i$th column of matrix $B$ has eight non-null coefficients corresponding to the velocity degrees of freedom located on the four edges of the element $E_i$. 

The upper part of system (\ref{global:linear:system}) permits to express the velocity vector ${U}_h$ in terms of the pressure $P_h$ in the following way
\begin{equation}\label{discrete:velocity}
{U}_h=-A^{-1}B\,P_h.
\end{equation}
Since $A$ is block-diagonal, this operation is locally performed on each block $A_i$. As a result, the velocity degrees of freedom associated to the $i$th vertex are expressed in terms of the pressure unknowns located at the centers of the elements sharing that vertex. Inserting (\ref{discrete:velocity}) into the lower part of (\ref{global:linear:system}), we get a cell-centered system for the pressures, i.e.,
\begin{equation}\label{pressure:system}
B^TA^{-1}B\,P_h=-F_h,
\end{equation}
whose matrix is symmetric and positive definite if we consider the symmetric quadrature rule (\ref{quadr:symm}). The discrete diffusion operator of MFMFE discretizations on logically rectangular grids has a 9-point stencil. In order to perform the LFA for the multigrid method presented in Section \ref{sec:2}, we need to derive the expression of the stencil coefficients for a homogeneous medium and a uniform Cartesian grid. Details are provided in the next subsection.

\subsection{The stencil coefficients for homogeneous media and Cartesian uniform grids}\label{sec:1.1}

Let us consider a homogeneous medium described by the permeability tensor
$$
K=\begin{bmatrix} a& c\\c& b\end{bmatrix},
$$
and a uniform Cartesian grid with mesh size $h$. In this case, $DF_{E}=hI$, where $I$ stands for the $2\times2$ identity matrix, and $J_E=h^2$. Hence, for every $E\in\mathcal{T}_h$,
$$\mathcal{K}_E=J_EDF_{E}^{-1}K(DF_{E}^{-1})^{T}=K.$$
The expression (\ref{block:vertex:matrix}) for the local matrix $A_i$ associated to the $i$th interior vertex is reduced to
$$A_i=\frac{1}{4(ab-c^2)}\begin{bmatrix}
2b\ &-c\ & 0\ &-c\ \\[0.5ex]-c & 2a & -c & 0\\[0.5ex] 0 & -c & 2b & -c\\[0.5ex]-c&0&-c&2a
\end{bmatrix}.
$$
If we insert this expression into (\ref{local:system:vertex}) and consider the local notations of Figure \ref{fig:int:node:2} (left), the volumetric fluxes associated to the $i$th interior vertex can be expressed in terms of pressure differences as
$$\begin{bmatrix}u_1\\[2ex] u_2\\[2ex] u_3 \\[2ex] u_4\end{bmatrix}=-\dfrac{1}{h}\begin{bmatrix}
a-\displaystyle\frac{c^2}{2b}\ &\displaystyle\frac{c}{2}\ & \dfrac{c^2}{2b}\ &\dfrac{c}{2}\ \\[2ex]\dfrac{c}{2} & b-\dfrac{c^2}{2a} & \dfrac{c}{2} & \dfrac{c^2}{2a}\\[2ex] \dfrac{c^2}{2b} & \dfrac{c}{2} & a-\dfrac{c^2}{2b} & \dfrac{c}{2}\\[2ex]\dfrac{c}{2}&\dfrac{c^2}{2a}&\dfrac{c}{2}&b-\dfrac{c^2}{2a}
\end{bmatrix}\begin{bmatrix}
p_2-p_1\\p_3-p_2\\p_3-p_4\\p_4-p_1\end{bmatrix}.$$

\begin{figure}[t]
	\begin{center}
		\unitlength=0.02cm
		\begin{picture}(200,200)
		
		\path(0,0)(180,0)(180,180)(0,180)(0,0)\path(60,0)(60,180)\path(120,0)(120,180) \path(0,60)(180,60)\path(0,120)(180,120)
		\put(25,26){$\times$}
		\put(85,26){$\times$}
		\put(145,26){$\times$}
		\put(25,85){$\times$}
		\put(85,85){$\times$}
		\put(145,85){$\times$}
		\put(25,145){$\times$}
		\put(85,145){$\times$}
		\put(145,145){$\times$}
		\put(18,16){\footnotesize{$p_{2}$}}
		\put(85,16){\footnotesize{$p_{3}$}}
		\put(152,16){\footnotesize{$p_{4}$}}
		\put(10,86){\footnotesize{$p_{9}$}}
		\put(85,75){\footnotesize{$p_{1}$}}
		\put(158,86){\footnotesize{$p_{5}$}}
		\put(15,158){\footnotesize{$p_{8}$}}
		\put(85,158){\footnotesize{$p_{7}$}}
		\put(154,158){\footnotesize{$p_{6}$}}
		\put(45,49){\footnotesize{$\mathbf{r}_1$}}
		\put(123,49){\footnotesize{$\mathbf{r}_2$}}
		\put(44,125){\footnotesize{$\mathbf{r}_4$}}
		\put(123,125){\footnotesize{$\mathbf{r}_3$}}
		\put(67,60){\circle*{5}}
		\put(111,60){\circle*{5}}
		\put(67,120){\circle*{5}}
		\put(111,120){\circle*{5}}
		\put(60,69){\circle*{5}}
		\put(120,69){\circle*{5}}
		\put(60,111){\circle*{5}}
		\put(120,111){\circle*{5}}
		\put(63,47){\footnotesize{$u_{1}$}}
		\put(103,47){\footnotesize{$u_{2}$}}
		\put(126,65){\footnotesize{$u_{3}$}}
		\put(126,107){\footnotesize{$u_{4}$}}
		\put(103,128){\footnotesize{$u_{5}$}}
		\put(63,128){\footnotesize{$u_{6}$}}
		\put(38,107){\footnotesize{$u_{7}$}}
		\put(38,65){\footnotesize{$u_{8}$}}
		\end{picture}\vspace*{0.5cm}
		\caption{Local notations for the 9-point stencil.}\vspace*{-0.5cm}
		\label{fig:stencil:cartesian}
	\end{center}
\end{figure}

Finally, we apply the expression (\ref{discrete:div}) to the central element of Figure \ref{fig:stencil:cartesian}. The local systems corresponding to the vertices $\mathbf{r}_1$, $\mathbf{r}_2$, $\mathbf{r}_3$ and $\mathbf{r}_4$ permit to express: $u_1$ and $u_8$ in terms of $p_1, p_2, p_3$ and $p_9$; $u_2$ and $u_3$ in terms of $p_1, p_3, p_4$ and $p_5$; $u_4$ and $u_5$ in terms of $p_1, p_5, p_6$ and $p_7$; and $u_6$ and $u_7$ in terms of $p_1, p_7, p_8$ and $p_9$. As a result, the 9-point stencil equation of system (\ref{pressure:system}) associated to the central pressure unknown $p_1$ is given by
$$\sum_{i=1}^9m_i\,p_i=\int_{E_1}f(\mathbf{x}),$$
where
\begin{align*}
&m_2=m_6=-\dfrac{c}{2}\left(1+\dfrac{c}{2a}+\dfrac{c}{2b}\right),& m_3=m_7=-b+\dfrac{c^2}{2}\left(\dfrac{1}{a}+\dfrac{1}{b}\right),\\[0ex] &m_4=m_8=\dfrac{c}{2}\left(1-\dfrac{c}{2a}-\dfrac{c}{2b}\right), & m_5=m_9=-a+\dfrac{c^2}{2}\left(\dfrac{1}{a}+\dfrac{1}{b}\right).
\end{align*}
The corresponding coefficient of $p_1$ is $m_1=-\textstyle\sum_{i=2}^9m_i$.

\section{Cell-centered multigrid and local Fourier analysis}
\label{sec:2}

In this section, we propose a blackbox cell-centered multigrid method for solving the MFMFE discretization of Darcy problem on logically rectangular grids. In addition, we apply a local Fourier analysis technique to study the convergence of the proposed solver.

\subsection{Cell-centered multigrid}\label{sec:2:mg}

The linear system \eqref{pressure:system} resulting from the MFMFE discretization described above demands efficient solvers for its solution. In this subsection, we design a cell-centered multigrid strategy \cite{MG2,WESSELING198885} that remains robust when applied to different logically rectangular grids with increasing levels of roughness. Remarkably, the proposed multigrid method acts as a blackbox, since the user simply provides a difference equation on a target grid and the code automatically generates the hierarchy of operators on coarser levels. In the sequel, we describe in detail the overall procedure.

Let $L_{\ell} P_{\ell} = F_{\ell}$ denote the linear system \eqref{pressure:system}, with $\ell$ representing the level of the target grid. In virtue of the blackbox multigrid idea, a hierarchy of smaller linear systems $$L_{k}P_k = F_k,\quad\ k=1,2,\ldots, \ell-1,$$  is built, where the corresponding linear operators on the coarser levels are constructed as
\begin{equation}\label{clc}
L_{k} = I_{k+1}^{k} L_{k+1} I_{k}^{k+1},
\end{equation}
$I_{k+1}^{k}$ and $I_{k}^{k+1}$ being appropriate restriction and prolongation operators, respectively. An iteration of a two-level cycle applies $\nu_1$ iterations of a relaxation procedure, followed by a coarse-level correction technique, and $\nu_2$ relaxation steps. In the coarse-level correction process, we first restrict the residual to the coarse level, then solve the coarse level problem exactly, and finally interpolate the obtained solution to the fine level in order to correct the previous approximation. This idea can be generalized to a multilevel approach by recursively applying the two-level algorithm.

The components of the proposed multigrid solver are explained next. The hierarchy of coarse level operators is constructed by Galerkin approximation as formulated in \eqref{clc}. As inter-grid transfer operators, we use a piecewise constant interpolation and the restriction operator by Wesseling/Khalil \cite{MG2}, that is,
\begin{equation}
\label{cp_wr}
I_{k}^{k+1} = 
\left]
\begin{array}{ccc}
1 &  & 1 \\ 
& \star & \\ 
1 & & 1
\end{array}
\right[_{k}^{k+1}, 
\qquad\quad 
I_{k+1}^{k}=\frac{1}{16}
\left[
\begin{array}{ccccc}
1 & 1 & & 0 & 0 \\ 
1 & 3 & & 2 & 0 \\ 
& & \star & & \\ 
0 & 2 & & 3 & 1 \\ 
0 & 0 & & 1 & 1
\end{array}
\right]_{k+1}^{k},
\end{equation}
where $\star$ denotes the position of a coarse level unknown. In the classical stencil notation for the prolongation, the numbers within the stencil represent the contribution of the coarse level unknown to the neighbouring fine level unknowns. In the stencil for the restriction, the numbers denote the weights corresponding to the fine level unknowns used to construct the coarse level value.

Regarding the relaxation procedure, a standard alternating line Gauss--Seidel is considered. One step of an elementary line relaxation updates all unknowns at the same \emph{grid line} simultaneously. In the case of structured rectangular grids, this implies considering the horizontal and vertical lines in the structured grid corresponding to the $x$- and $y$-line relaxation, respectively. However, in the case of logically rectangular grids, this has to be defined more precisely. Since the index set $\{(i,j),\,i=1,2,\ldots,n,\,j=1,2,\ldots,m\}$ of a logically rectangular grid has the same structure as that of a rectangular grid, we can define the ``$x$-line'' smoother as the relaxation that simultaneously updates all the unknowns located at grid points with a fixed $i$-index. In turn, the ``$y$-line'' smoother updates together those unknowns located at grid points with a fixed $j$-index. 
Finally, the alternating version of the line relaxation step consists of an $x$-line relaxation iteration followed by a $y$-line relaxation step.

\begin{remark}
Since our approach applies geometric multigrid on logically rectangular grids, an easy parallelization of the algorithm can be obtained by grid partitioning \cite{TOS01}. The domain can be split into several blocks and each block is then sent to a processor of a parallel computer. 
\end{remark}

\subsection{Local Fourier analysis} \label{sec:2_LFA}

In order to apply the LFA \cite{Bra77,Bra94}, we restrict ourselves to studying the Darcy system in a homogenous medium, discretized by an MFMFE scheme on a Cartesian uniform grid, whose mesh size is $h$ in both directions. In this case, we obtain a discrete diffusion operator with constant coefficients, which can be represented by a constant stencil that is common to all the interior points. This stencil was derived in Section \ref{sec:1}. Such an operator is assumed to be defined on an infinite grid $G_h$, thus neglecting the effect of the boundary conditions.

Based on these assumptions, all the operators involved in the multigrid method are extended to the infinite grid, and the error grid function can be written as a formal linear combination of the so-called Fourier modes, i.e., $$\varphi_h({\boldsymbol \theta}, {\mathbf x}) = e^{\imath {\boldsymbol \theta}{\mathbf x}/{\mathbf h}} = e^{\imath \theta_1x_1/h} e^{\imath \theta_2 x_2/ h}, \ \hbox{where } \imath=\sqrt{-1},$$ which span the Fourier space $${\mathcal F}(G_h) = \{ \varphi_h({\boldsymbol \theta}, {\mathbf x}) \; | \; {\boldsymbol \theta} = (\theta_1, \theta_2) \in \Theta_h= (-\pi,\pi]^2, \; {\mathbf x} = (x_1,x_2) \in G_h\}.$$ 
The main idea of the LFA is to study how the two-grid operator acts on this expression of the error based on the Fourier components. Since the discrete diffusion operator satisfies the LFA assumptions, the Fourier modes are its eigenfunctions. More precisely, we obtain 
$$L_{\ell} \varphi_h({\boldsymbol \theta}, {\mathbf x}) = \widetilde{L}_{\ell}({\boldsymbol \theta}) \varphi_h({\boldsymbol \theta}, {\mathbf x}),$$
where $\widetilde{L}_{\ell}({\boldsymbol \theta})$ is called the \textit{symbol} of $L_{\ell}$.
In our case, recalling the stencil of the discrete operator obtained at the end of Section \ref{sec:1}, we have 
$$\widetilde{L}_{\ell}({\boldsymbol \theta}) = m_1 + 2m_2\cos(\theta_1+\theta_2) + 2m_3\cos(\theta_2) + 2m_4\cos(\theta_1-\theta_2)+ 2m_5\cos(\theta_1).$$
It is well known that the same holds true for standard relaxation operators as those considered here, so that we can compute $\widetilde{S}_{\ell}({\boldsymbol \theta})$ and, consequently, the LFA smoothing factor is given by
$$\mu = \sup_{\Theta_h\backslash \Theta_h^{\mathrm{low}}} \rho(\widetilde{S}_{\ell}({\boldsymbol \theta})),$$
where $\Theta_h^{\mathrm{low}} = (-\pi/2,\pi/2]^2$ is the set of low frequencies associated with standard coarsening. However, the inter-grid transfer operators couple Fourier modes, thus leaving invariant the subspaces of $2h$-harmonics,
$${\mathcal F}^2({\boldsymbol \theta}^{00}) = span\left\{\varphi_h({\boldsymbol \theta}^{00},\cdot),\varphi_h({\boldsymbol \theta}^{11},\cdot),\varphi_h({\boldsymbol \theta}^{10},\cdot),\varphi_h({\boldsymbol \theta}^{01},\cdot)\right\},$$
where ${\boldsymbol \theta}^{\alpha\beta} = {\boldsymbol \theta}^{00} - (\alpha \, \hbox{sign} (\theta_1), \beta \, \hbox{sign} (\theta_2))\pi$, for $\alpha,\,\beta\in\{0,1\}$. In this way, the Fourier symbols of the restriction and prolongation operators are matrices of sizes $1\times 4$ and $4\times 1$, respectively, given by 
\begin{eqnarray*}
	\widetilde{I}_{\ell}^{\ell-1} ({\boldsymbol \theta}) &=& \left[\widetilde{R}({\boldsymbol \theta}^{00}) \; \widetilde{R}({\boldsymbol \theta}^{11}) \; \widetilde{R}({\boldsymbol \theta}^{10})\; \widetilde{R}({\boldsymbol \theta}^{01})\right],\\
	\widetilde{I}_{\ell-1}^{\ell} ({\boldsymbol \theta}) &=& \left[\widetilde{P}({\boldsymbol \theta}^{00}) \; \widetilde{P}({\boldsymbol \theta}^{11}) \; \widetilde{P}({\boldsymbol \theta}^{10})\; \widetilde{P}({\boldsymbol \theta}^{01})\right]^T, 
\end{eqnarray*}
where 
\begin{eqnarray*}
	\widetilde{R}({\boldsymbol \theta}^{\alpha\beta}) &=& \cos(\theta_1^{\alpha \beta}/2) \, \cos(\theta_2^{\alpha \beta}/2)\, \cos((\theta_1^{\alpha \beta}-\theta_2^{\alpha \beta})/2),\\
	\widetilde{P}({\boldsymbol \theta}^{\alpha\beta}) &=& \cos(\theta_1^{\alpha \beta}/2) \, \cos(\theta_2^{\alpha \beta}/2),
\end{eqnarray*} 
with $\alpha,\,\beta\in\{0,1\}$. Using the preceding symbols and taking into account that the symbol of the Galerkin operator is given by $$\widetilde{L}_{\ell-1}({\boldsymbol \theta}) = \widetilde{I}_{\ell}^{\ell-1} ({\boldsymbol \theta}) \widetilde{L}_{\ell}({\boldsymbol \theta}) \widetilde{I}^{\ell}_{\ell-1} ({\boldsymbol \theta}),$$ we can construct the Fourier representation of the two-level operator $M_{\ell}^{\ell-1}$ as follows
$$\widetilde{M}_{\ell}^{\ell-1}({\boldsymbol \theta}) = \widetilde{S}_{\ell}^{\nu_2}({\boldsymbol \theta}) \left(\widetilde{I}_{\ell}({\boldsymbol \theta}) - \widetilde{I}_{\ell-1}^{\ell} ({\boldsymbol \theta}) \widetilde{L}_{\ell-1}^{-1}({\boldsymbol \theta}) \widetilde{I}_{\ell}^{\ell-1} ({\boldsymbol \theta}) \widetilde{L}_{\ell} ({\boldsymbol \theta}) \right) \widetilde{S}_{\ell}^{\nu_1}({\boldsymbol \theta}),$$
where $\widetilde{I}_{\ell}({\boldsymbol \theta})$ denotes the $4\times 4$ identity matrix.
Finally, we can estimate the spectral radius of the two-level operator, which provides a prediction of the asymptotic convergence factor of the method, as 
$$\rho_{2g} = \sup_{{\boldsymbol \theta}\in \Theta_h^{\mathrm{low}}} \rho \left(\widetilde{M}_{\ell}^{\ell-1}({\boldsymbol \theta})\right).$$

The LFA described above allows us to estimate the convergence rates of the proposed algorithm very accurately. The main results of the analysis are shown in Table \ref{lfa_results}, which contains the smoothing factors, $\mu$, and two-grid convergence factors, $\rho_{2g}$, predicted by the LFA. The table further displays the real asymptotic convergence factors, $\rho_h$, obtained by using $W$-cycles. The test considers one smoothing step of two different smoothers, Gauss--Seidel and alternating line relaxations, and three different tensor samples. In particular, we take $K_1$ to be the identity matrix, and $K_2$ and $K_3$ the following full tensors
$$
K_2=\begin{bmatrix}
4 & 1 \\
1 & 4
\end{bmatrix},\qquad K_3=\begin{bmatrix}
2 & 1 \\
1 & 10000
\end{bmatrix}.
$$
Notice that $K_2$ is isotropic, whereas $K_3$ is highly anisotropic. We can observe in Table \ref{lfa_results} that the two-level convergence factors provided by the LFA match the experimentally computed asymptotic convergence factors very accurately. It can also be observed that, in the case of the highly anisotropic tensor, the multigrid method based on the Gauss--Seidel smoother does not converge, and the alternating line relaxation is needed to obtain a successful result.  These results are in accordance with previous findings that state that block-wise smoothers are much more efficient than point-wise smoothers when considering anisotropic problems \cite{TOS01}. Based on these observations, we will choose the alternating line smoother to obtain a robust blackbox-type multigrid solver for any arbitrary tensor $K$.
\begin{table}[t]
	\caption{LFA-predicted smoothing factors ($\mu$) and two-grid convergence factors ($\rho_{2g}$), together with the experimentally obtained asymptotic convergence factors ($\rho_h$), for three different tensors.}
	\label{lfa_results}
	\begin{center}
		\begin{tabular}{ccccccc}\cline{2-7}
			& \multicolumn{3}{c}{Gauss--Seidel}&\multicolumn{3}{c}{Alternating line}\\
			\cline{2-7}
			& $\mu$ & $\rho_{2g}$ & $\rho_h$ & $\mu$ & $\rho_{2g}$ & $\rho_h$\\
			\hline
			$K_1$ & 0.50 & 0.40 & 0.40 & 0.15 & 0.11 & 0.12\\
			$K_2$ & 0.50 & 0.42 & 0.41 & 0.15 & 0.19 & 0.18\\
			$K_3$ & 1.00 & 1.00 & -- & 0.45 & 0.22 & 0.17\\ 
			\hline
		\end{tabular}
	\end{center}
\end{table}

\section{Numerical experiments}
\label{sec:3}

\subsection{A test with a known analytical solution}
\label{subsec:3_1}

Let us consider a boundary value problem of type (\ref{cont:problem}), where $\Omega=(0,1)^2$, $\partial\Omega=\Gamma_D$ and tensor $K$ is given by
$$
K=\begin{bmatrix}
5 & 3 \\
3 & 7
\end{bmatrix}.
$$
Data functions $f$ and $g$ are defined in such a way that the exact solution of the problem is $p(x,y)=\sin(\pi x)^2\sin(2\pi y)$. 

\begin{figure}[t]\vspace*{0.36cm}
	\begin{center}
		\begin{minipage}[t]{0.44\textwidth}
			\begin{center}\includegraphics[scale=0.16]{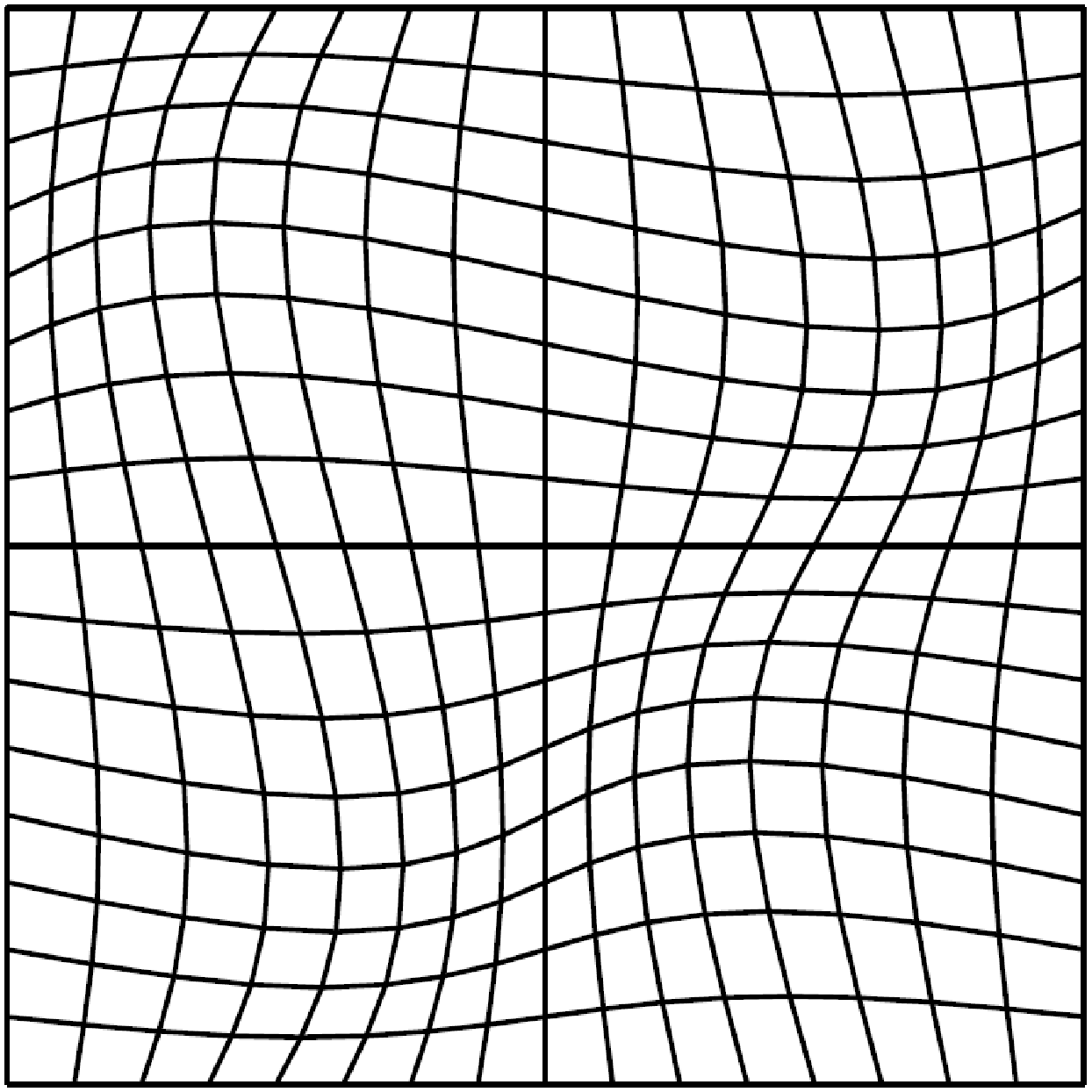}\\[1ex]{\footnotesize (a) Smooth mesh}\end{center}
		\end{minipage}
		\hspace*{-0.92cm}
		\begin{minipage}[t]{0.44\textwidth}\begin{center}\vspace*{-3.7cm}\includegraphics[scale=0.37]{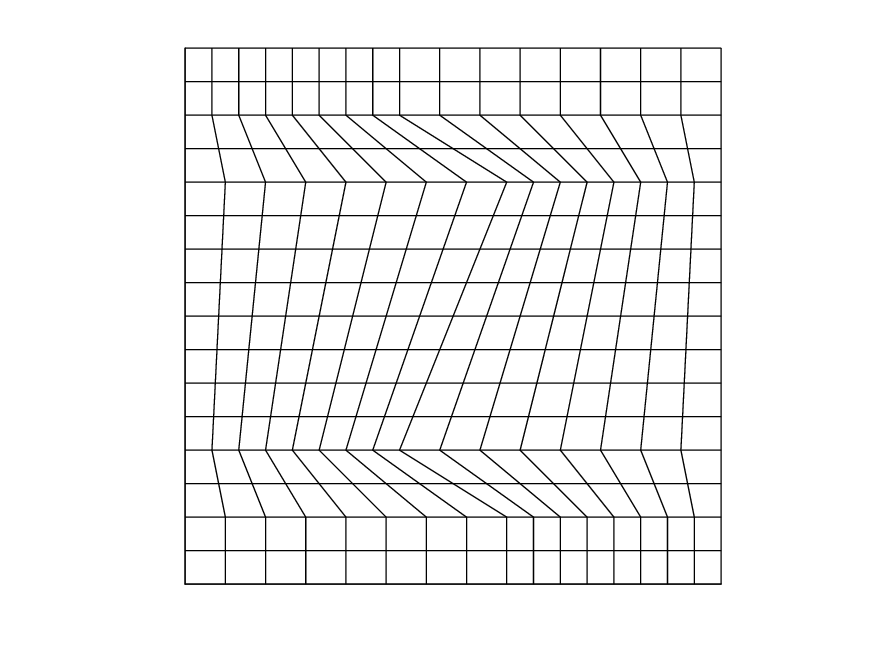}\\[-2.5ex]{\footnotesize (b) Kershaw mesh}\end{center}
		\end{minipage}\\[2ex]
		\begin{minipage}[t]{0.44\textwidth}\begin{center}\includegraphics[scale=0.215]{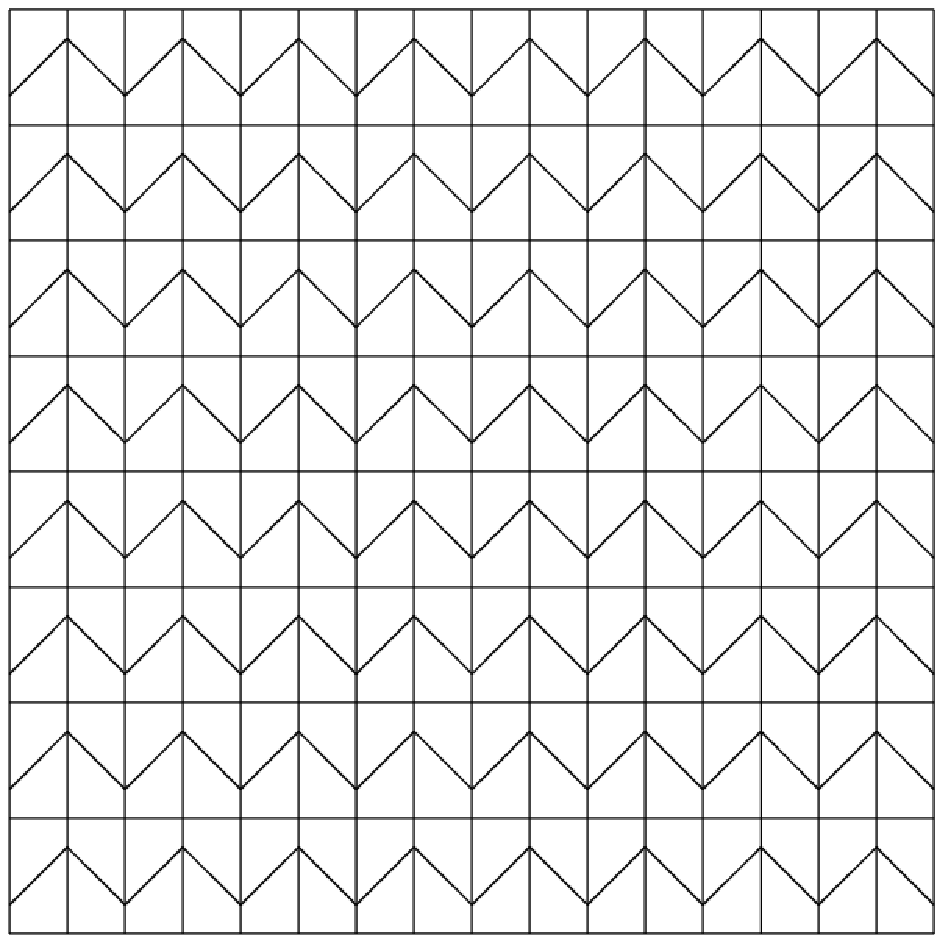}\\[1ex]{\footnotesize (c) $h$-perturbed mesh}\end{center}
		\end{minipage}
		\hspace*{-0.8cm}
		\begin{minipage}[t]{0.44\textwidth}\begin{center}\includegraphics[scale=0.16]{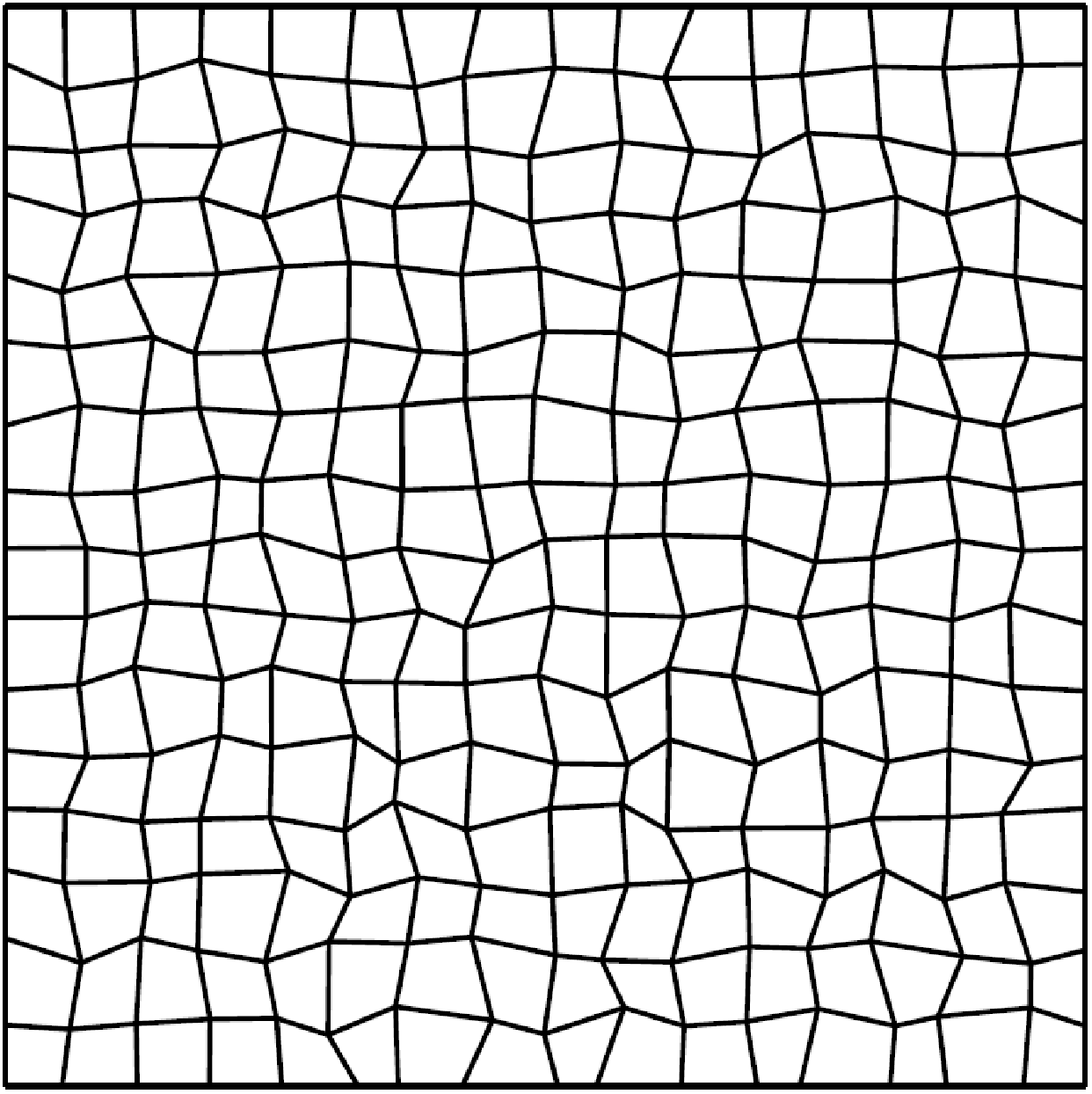}\\[1ex]{\footnotesize (d) Randomly $h$-perturbed mesh}\end{center}
		\end{minipage}\vspace*{0.2cm}
		\caption{Some logically rectangular grids considered in the numerical examples.}\label{fig:meshes}
	\end{center}
\end{figure}

The spatial domain is discretized using four different logically rectangular meshes that contain $N\times N$ elements. The first one is a family of smooth meshes defined to be a $C^{\infty}$-mapping of successively refined uniform meshes on the unit square, i.e.,
\begin{align*}
x_{i,j}&=\hat{x}_{i,j}+\frac{3}{50}\sin(2\pi\hat{x}_{i,j})\sin(2\pi\hat{y}_{i,j}),\\
y_{i,j}&=\hat{y}_{i,j}-\frac{1}{20}\sin(2\pi\hat{x}_{i,j})\sin(2\pi\hat{y}_{i,j}),
\end{align*}
where $\hat{x}_{i,j}$ and $\hat{y}_{i,j}$ denote the spatial coordinates of the vertices on the uniform mesh, and $x_{i,j}$ and $y_{i,j}$ are their counterparts on the smooth mesh, for $i=1,2,\ldots,N+1$ and $j=1,2,\ldots,N+1$. An illustration of this type of meshes is given in Figure \ref{fig:meshes}(a). Next, we consider a set of Kershaw-type meshes \cite{ker:81,sha:ste:96}, which contain certain highly skewed zones that are displayed on Figure \ref{fig:meshes}(b). The $h$-perturbed grid represented in Figure \ref{fig:meshes}(c) belongs to a family of meshes described in detail in \cite{arn:bof:fal:05}. Finally, we consider a family of randomly $h$-perturbed meshes consisting of highly distorted quadrilaterals. As shown in Figure \ref{fig:meshes}(d), each of these meshes is generated by perturbing the vertices of a uniform mesh by a distance of size $O(h)$ in a random direction. More specifically, the spatial coordinates of the vertices of the randomly $h$-perturbed mesh can be obtained as
\begin{align*}
x_{i,j}&=\hat{x}_{i,j}-\frac{13}{10}h+\frac{3}{5}hr_x^{i,j},\\
y_{i,j}&=\hat{y}_{i,j}-\frac{13}{10}h+\frac{3}{5}hr_y^{i,j},
\end{align*}
where $h=1/N$, and $r_x^{i,j}$ and $r_y^{i,j}$ are pseudo-random numbers uniformly distributed on the interval $(0,1)$. It is interesting to note that, from the four meshes under consideration, just the smooth and the Kershaw-type grids are composed of $\mathcal{O}(h^2)$-parallelograms.

Convergence results for the MFMFE discretizations on different types of meshes are extensively reported in the literature (e.g., \cite{aav:eig:kla:whe:yot:07,whe:xue:yot:12a,whe:yot:06}). For the sake of illustration, we include here some convergence tests for both the symmetric and non-symmetric MFMFE schemes on the Kershaw and $h$-perturbed grids. As reported below, the symmetric method shows a good convergence on the Kersahw mesh, while the non-symmetric variant is required when considering the $h$-perturbed grid. In this setting, we compute the global errors for the pressure and velocity variables using different norms, i.e.,
\begin{align*}
E^p_{h}&=\|p-p_h\|_{L^2},&
\widehat{E}^p_{h}&=\|r_hp-P_h\|_{\ell^2},\\
{E}^{\mathbf{u}}_{h}&=\|\Pi_h \mathbf{u}-\mathbf{u}_h\|_{L^2},&
\widehat{E}^{\mathbf{u}}_{h}&=\|\mathbf{u}-\mathbf{u}_h\|_{\mathcal{F}_h},
\end{align*}
where $r_h$ denotes the restriction operator to the center of the cells, and $\Pi_h$ denotes the projection operator onto the space $V_h$ \cite{whe:yot:06}. We denote by $\|\cdot\|_{L^2}$ the corresponding norm in either $L^2(\Omega)$ or $(L^2(\Omega))^2$, and by $\|\cdot\|_{\ell^2}$ the discrete $L^2$-norm in $\mathcal{H}_p$. Further, we define the edge-based norm $\|\cdot\|_{\mathcal{F}_h}$ as
\begin{equation*}
\| \textbf{v} \|^{2}_{\mathcal{F}_{h}} = \sum_{E\in \mathcal{T}_{h}}\sum_{e\in \partial E} \frac{|E|}{|e|}\| \textbf{v}\cdot\textbf{n}_{e} \|^{2}_{e},
\end{equation*}
where $|E|$ and $|e|$ refer to the area of element $E$ and the length of edge $e$, respectively. The integrals involved in the pressure errors $E^p_{h}$ are approximated element-wise by a 9-point Gaussian quadrature formula. In turn, the integrals arising in the velocity error ${E}^{\mathbf{u}}_{h}$ are approximated by the trapezoidal quadrature rule, while those involved in the face error $\widehat{E}^{\mathbf{u}}_{h}$ are computed by a high-order Gaussian quadrature rule.

Table \ref{table:kershaw} shows the global errors and numerical orders of convergence of the method when applied to the family of Kershaw-type meshes. As predicted by the theory \cite{whe:yot:06}, we observe first-order convergence for the pressure and the velocity errors in the $L^2$-norm, as well as for the velocity errors on the element edges. Moreover, we obtain second-order superconvergence for the pressure at the cell centers.
\begin{table}[t]
	\caption{Global errors and numerical orders of convergence for the symmetric MFMFE method on the Kershaw meshes ($h_0=2^{-5}$).}
	\label{table:kershaw}
	\renewcommand{\arraystretch}{1.2}
	\begin{tabular}{l|cccccccc} 
		\hline\\[-3ex]
		$h$ & $E^p_{h}$ & Rate & $\hat{E}^p_{h}$ & Rate & $E^{\mathbf{u}}_{h}$ & Rate & $\hat{E}^{\mathbf{u}}_{h}$ & Rate \\
		\hhline{|---------|}
		$h_0$ & 2.959e-02 & -- & 1.857e-03 & -- & 1.103e{\small+}00& -- & 9.141e-01& -- \\
		$h_0/2$& 1.479e-02 & 1.000  & 4.739e-04 & 1.970 & 5.519e-01& 0.999 & 4.524e-01& 1.014 \\
		$h_0/2^2$& 7.396e-03 & 1.000 & 1.192e-04 & 1.991 & 2.763e-01& 0.998 & 2.258e-01& 1.002 \\
		$h_0/2^3$& 3.698e-03 & 1.000 & 2.986e-05 & 1.997 & 1.383e-01& 0.999 & 1.129e-01& 1.000 \\
		$h_0/2^4$& 1.851e-03 & 0.998 & 7.470e-06 & 1.999 & 6.916e-02& 1.000 &  5.633e-02& 1.003 \\
		\hline
	\end{tabular}
\end{table}

Table \ref{table:boffi:sym} shows the global errors and the numerical orders of convergence on the family of $h$-perturbed grids depicted in Figure \ref{fig:meshes}(c). The numerical results show that the convergence of the symmetric MFMFE method deteriorates on this type of grids.
	\begin{table}[t]
			\caption{Global errors and numerical orders of convergence for the symmetric MFMFE method on the $h$-perturbed meshes ($h_0=2^{-5}$).}
			\label{table:boffi:sym}
			\renewcommand{\arraystretch}{1.2}
			\begin{tabular}{l|cccccccc}
				\hline\\[-3ex]
				$h$ & $E^p_{h}$ & Rate & $\hat{E}^p_{h}$ & Rate & $E^{\mathbf{u}}_{h}$ & Rate & $\hat{E}^{\mathbf{u}}_{h}$ & Rate \\
				\hhline{|---------|}
				$h_0$ & 3.400e-02 & -- & 7.831e-03 & -- & 2.712e{\small+}00& -- & 2.349e{\small +}00& -- \\
				$h_0/2$& 1.702e-02 & 0.998 & 3.821e-03 & 1.035 & 1.804e{\small +}00& 0.588 & 1.565e{\small +}00& 0.586  \\
				$h_0/2^2$& 8.513e-03 & 0.999 & 1.893e-03 & 1.013 & 1.234e{\small +}00& 0.548 & 1.072e{\small +}00& 0.546 \\
				$h_0/2^3$& 4.258e-03 & 0.999 & 9.432e-04 & 1.005 & 8.575e-01& 0.525 & 7.452e-01& 0.525\\
				$h_0/2^4$& 2.129e-03 & 1.000 &  4.708e-04 & 1.002 & 6.010e-01& 0.513 &  5.224e-01& 0.512\\
				\hline
		\end{tabular}
	\end{table}
	In turn, Table \ref{table:boffi:non:sym} shows the numerical results obtained when considering the non-sym\-metric MFMFE method derived by using the non-symmetric quadrature rule (\ref{quadr:non:symm}). In accordance with \cite{whe:xue:yot:12a}, we observe first-order convergence for both the velocity and the pressure errors in the $L^2$-norm, as well as for the velocity errors on the element edges. Finally, second-order superconvergence is obtained for the pressure errors at the cell centers.
\begin{table}[t]
		\caption{Global errors and numerical orders of convergence for the non-symmetric MFMFE method on the $h$-perturbed meshes ($h_0=2^{-5}$).}
		\label{table:boffi:non:sym}
		\renewcommand{\arraystretch}{1.2}
		\begin{tabular}{l|cccccccc}
			\hline\\[-3ex]
			$h$ & $E^p_{h}$ & Rate & $\hat{E}^p_{h}$ & Rate & $E^{\mathbf{u}}_{h}$ & Rate & $\hat{E}^{\mathbf{u}}_{h}$ & Rate \\
			\hhline{|---------|}
			$h_0$ & 3.325e-02 & -- & 1.658e-03 & -- & 1.370e{\small+}00& -- & 1.072e{\small+}00& -- \\
			$h_0/2$& 1.662e-02 & 1.000  & 4.163e-04 & 1.994 & 6.855e-01& 0.999 & 5.339e-01& 1.006 \\
			$h_0/2^2$& 8.308e-03 & 1.000 & 1.042e-04 & 1.998 & 3.428e-01& 1.000 & 2.667e-01& 1.001 \\
			$h_0/2^3$& 4.154e-03 & 1.000 & 2.605e-05 & 2.000 & 1.714e-01& 1.000 & 1.333e-01& 1.000 \\
			$h_0/2^4$& 2.077e-03 & 1.000 &  6.513e-06 & 1.999 & 8.571e-02& 1.000 &  6.665e-02& 1.000\\
			\hline
	\end{tabular}
\end{table}

Next, let us show the performance of the multigrid method when considering different types of cycles ($V$-, $F$- and $W$-cycles), different families of meshes and different grid sizes. In Table \ref{table:iterations}, we show the number of multigrid iterations required to obtain a residual smaller than $10^{-9}$. In the case of considering random grids, the last three columns of Table \ref{table:iterations} show the rounding to the closest natural number of the average number of multigrid iterations required when performing 100 realizations of the random mesh. It is worth pointing out that in the case of considering Kershaw meshes, we propose to use a relaxation parameter $w=0.6$ in the smoothing procedure. The multigrid solver is shown to be robust in all the cases with respect to the spatial discretization parameter $h$. Moreover, it is shown to be very efficient, since the number of iterations required to satisfy the stopping criterion  in all the cases is not very high.  It is also interesting to notice that both $F$- and $W$-cycles show a better convergence behaviour than $V$-cycles on rough grids such as Kershaw meshes. Finally, since the computational cost of $F$-cycles is smaller than that of $W$-cycles, from now on we shall consider $F$-cycles for the multigrid solver.
\begin{table}[t]
		\caption{Number of iterations of the multigrid method when considering different types of cycles and various families of logically rectangular meshes ($h_0=2^{-5}$) in the test with known analytical solution.}
		\label{table:iterations}
		\renewcommand{\arraystretch}{1.2}
		\begin{center}
			\begin{tabular}{c|ccc|ccc|ccc|ccc}
				\hline\\[-3ex]
				Mesh&\multicolumn{3}{c|}{Smooth}&\multicolumn{3}{c|}{Kershaw}&\multicolumn{3}{c|}{$h$-perturbed}& \multicolumn{3}{c}{Randomly $h$-perturbed}  \\\hline
				Cycle&$V$& $F$& $W$& $V$& $F$& $W$& $V$& $F$& $W$& $V$& $F$& $W$
				\\\hline
				$h_0$&10 & 9 & 9 & 9 & 9 & 9 & 8 & 8 & 8 & \hspace*{0.37cm}9\hspace*{0.37cm} &  \hspace*{0.15cm}8\hspace*{0.15cm}& 8
				\\
				$h_0/2$&11 & 9 & 9 & 11 & 10 & 10 & 8 & 7 & 7 & 9 & 8 & 8 
				\\
				$h_0/2^2$& 11 & 7 & 7 & 14 & 13 & 13 & 8 & 6 & 6 & 9  & 8 & 8 
				\\
				$h_0/2^3$& 11 & 5 & 5 & 16 & 13 & 13 & 8 & 5 & 5 & 9  & 7 & 7  
				\\
				$h_0/2^4$& 10 & 4 & 4 & 17 & 13 & 13&8 &5 &5 & 9 & 7 & 7 
				\\\hline
				
			\end{tabular}
		\end{center}
\end{table}

\subsection{Several tests with discontinuous permeability tensors}
\label{subsec:3_2}

Let us next consider a flow problem through a system that contains an impermeable streak. Different variants of this porous media example have been considered as test problems for various spatial discretizations in \cite{arr:por:yot:14,dur:93,dur:94,hym:sha:ste:97}. In our case, this numerical example permits us to test the behaviour of the new multigrid method in the solution of a problem that contains irregularly shaped strata and abrupt variations in permeability.

In particular, we consider problem (\ref{cont:problem}) with $\Omega=(0,1)^2$, Dirichlet boundary conditions with $g=1-x$ on $\Gamma_D=\partial\Omega$, and $f=0$. The flow domain contains a low-permeability region which is delimited by two curves. The top curve is chosen to be an arc of a circle with center at $(0.1,-0.4)$ and radius equal to $1.2$, while the bottom curve is an arc of a circle with the same center and radius equal to $1.1$. The permeability throughout the domain is uniform and isotropic ($K=I$), except in the low-permeability streak. In this region, it is such that the parallel component to the local streak orientation is equal to $0.1$, and the normal component to the local streak orientation is equal to 0.001. For more details, see \cite{arr:por:yot:14}.

\begin{figure}[t]
	\begin{center}
		\hspace*{-0.9cm}\begin{minipage}[t]{0.3\textwidth}
			\begin{center}
				\includegraphics[scale=0.32]{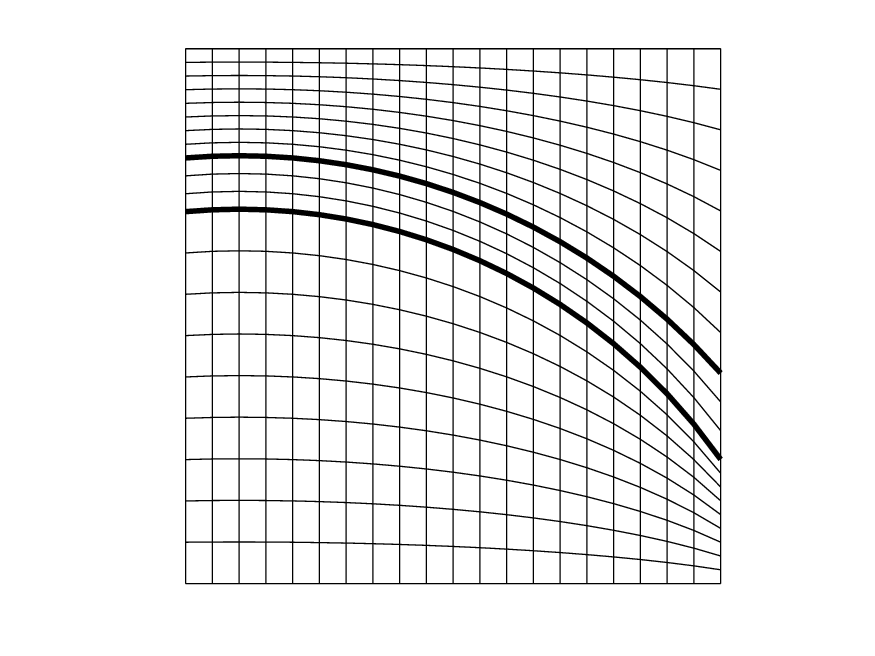}
			\end{center}
		\end{minipage}
		\hspace*{0.3cm}
		\begin{minipage}[t]{0.3\textwidth}
			\begin{center}
				\includegraphics[scale=0.32]{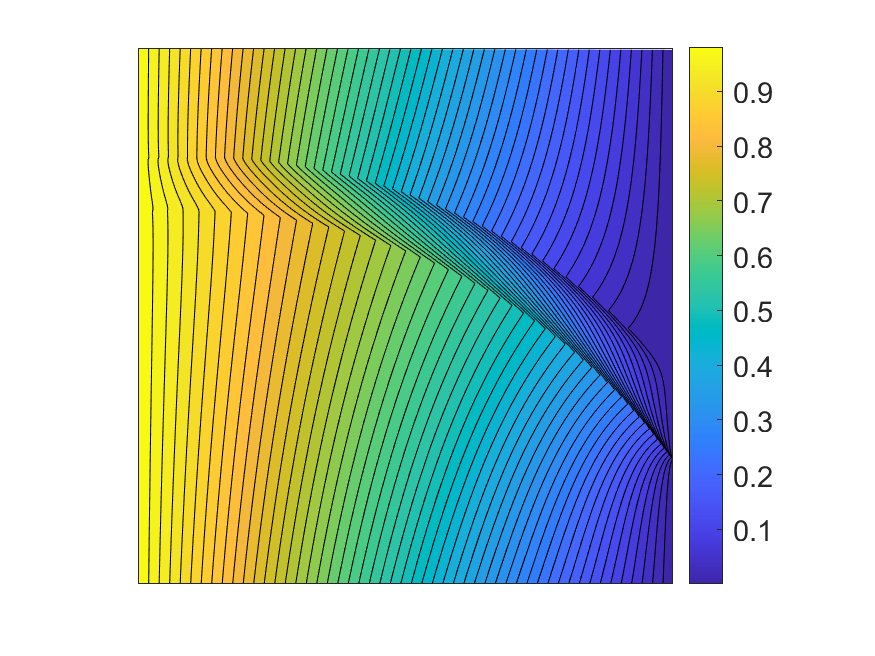}
			\end{center}
		\end{minipage}
		\hspace*{0.7cm}
		\begin{minipage}[t]{0.3\textwidth}
			\begin{center}\vspace*{-0.187cm}
				\includegraphics[scale=0.355]{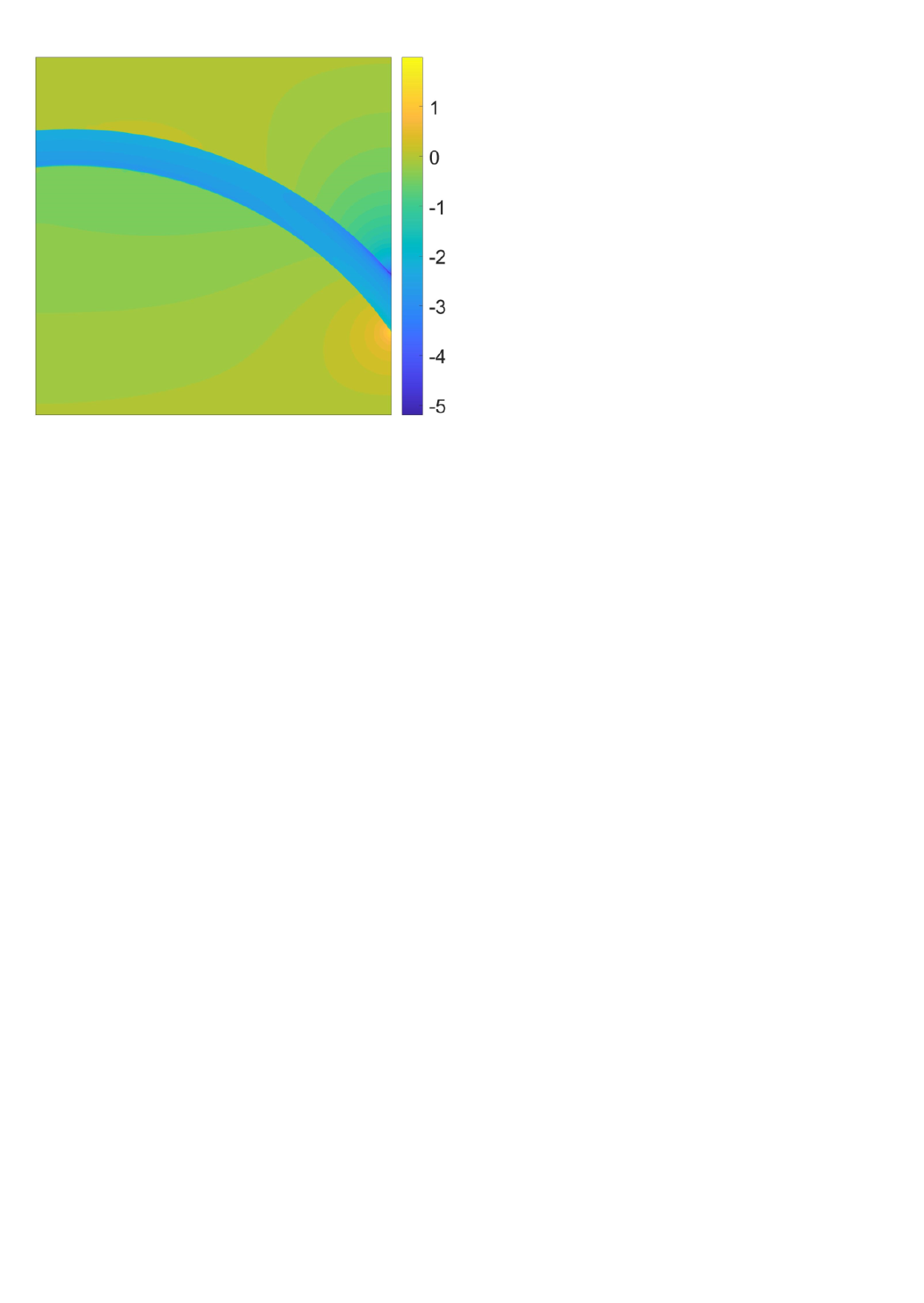}
			\end{center}
		\end{minipage}\vspace*{-6.8cm}
	\end{center}
\caption{Logically rectangular mesh (left), numerical pressure (center) and  logarithm of  the  norm  of  the  numerical  velocity  (right) for the flow problem that contains an impermeable streak.}\label{fig:porous:media}
\end{figure}
		
\begin{figure}[t]
	\begin{center}
				\includegraphics[scale=0.32]{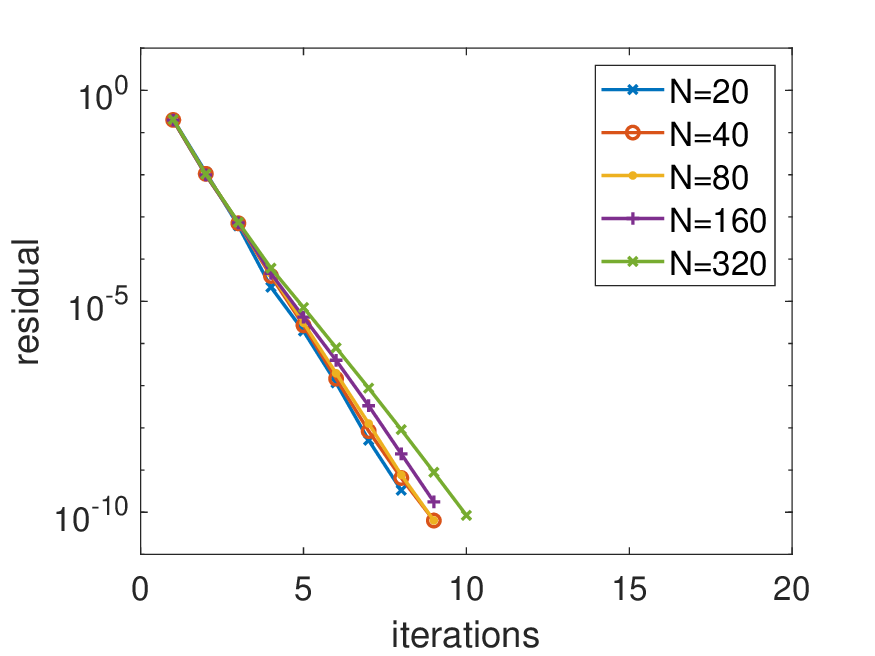}
			\end{center}
	\caption{ Convergence history of the multigrid method for the flow problem that contains an impermeable streak.}\label{fig:porous:media:conv:his}
\end{figure}

Figure \ref{fig:porous:media} (left) shows the geometry of the flow domain and the logically rectangular grid used in the discretization when $N=20$. Observe that the grid is adapted to the geometry of the low-permeability streak, depicted in the figure by bold curves.
Figure \ref{fig:porous:media} (center) displays the numerical pressure and Figure \ref{fig:porous:media} (right) shows the logarithm of the norm of the numerical velocity field when considering $N=320$. As expected from the physical configuration, no flow enters the streak, which is in accordance with the numerical results obtained in the aforementioned works.
Regarding the performance of the multigrid solver, Figure \ref{fig:porous:media:conv:his} shows the residual versus the number of iterations until the initial residual is reduced by a factor of $10^{-10}$. In this case, where we consider a relaxation parameter $w=0.6$ in the smoothing procedure, the multigrid solver is also shown to be robust with respect to the spatial discretization parameter $h$.  

\begin{figure}[t]
		\begin{center}
			\hspace*{-0.7cm}\begin{minipage}[t]{0.3\textwidth}
				\begin{center}
					\includegraphics[scale=0.32]{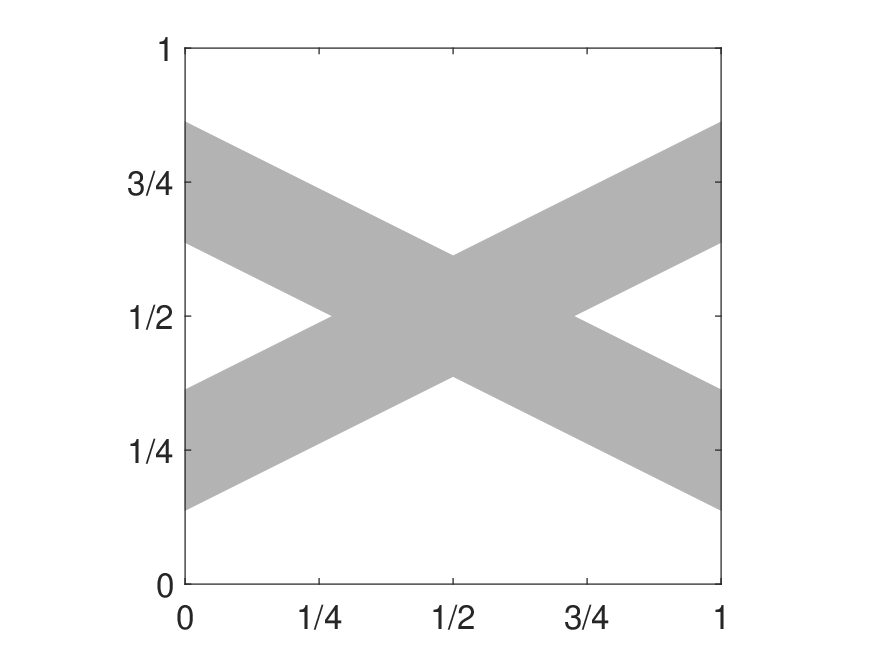}
				\end{center}
			\end{minipage}
			\hspace*{0.3cm}
			\begin{minipage}[t]{0.3\textwidth}
				\begin{center}
					\includegraphics[scale=0.32]{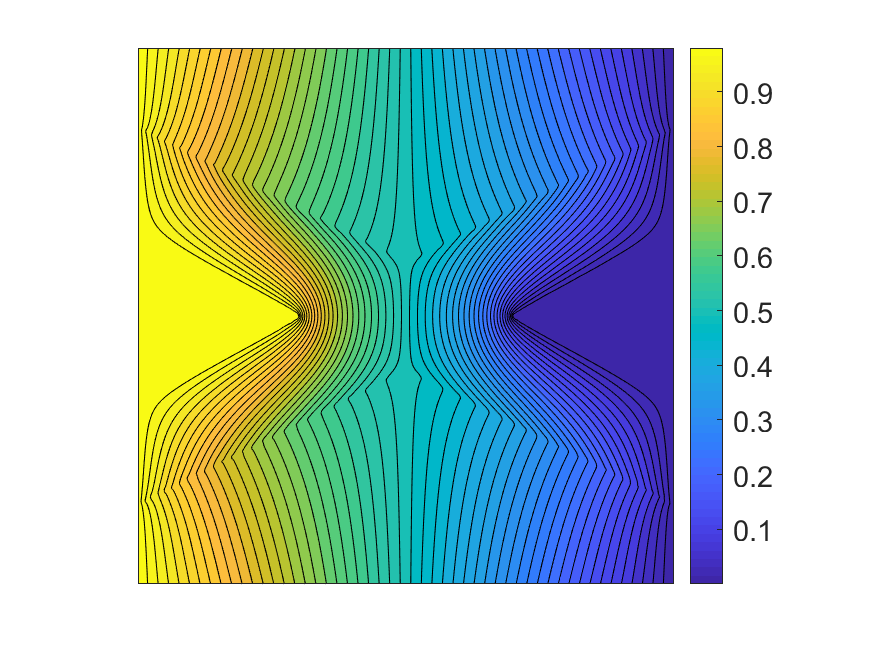}
				\end{center}
			\end{minipage}
			\hspace*{0.7cm}
			\begin{minipage}[t]{0.3\textwidth}
				\begin{center}\vspace*{-0.187cm}
					\includegraphics[scale=0.355]{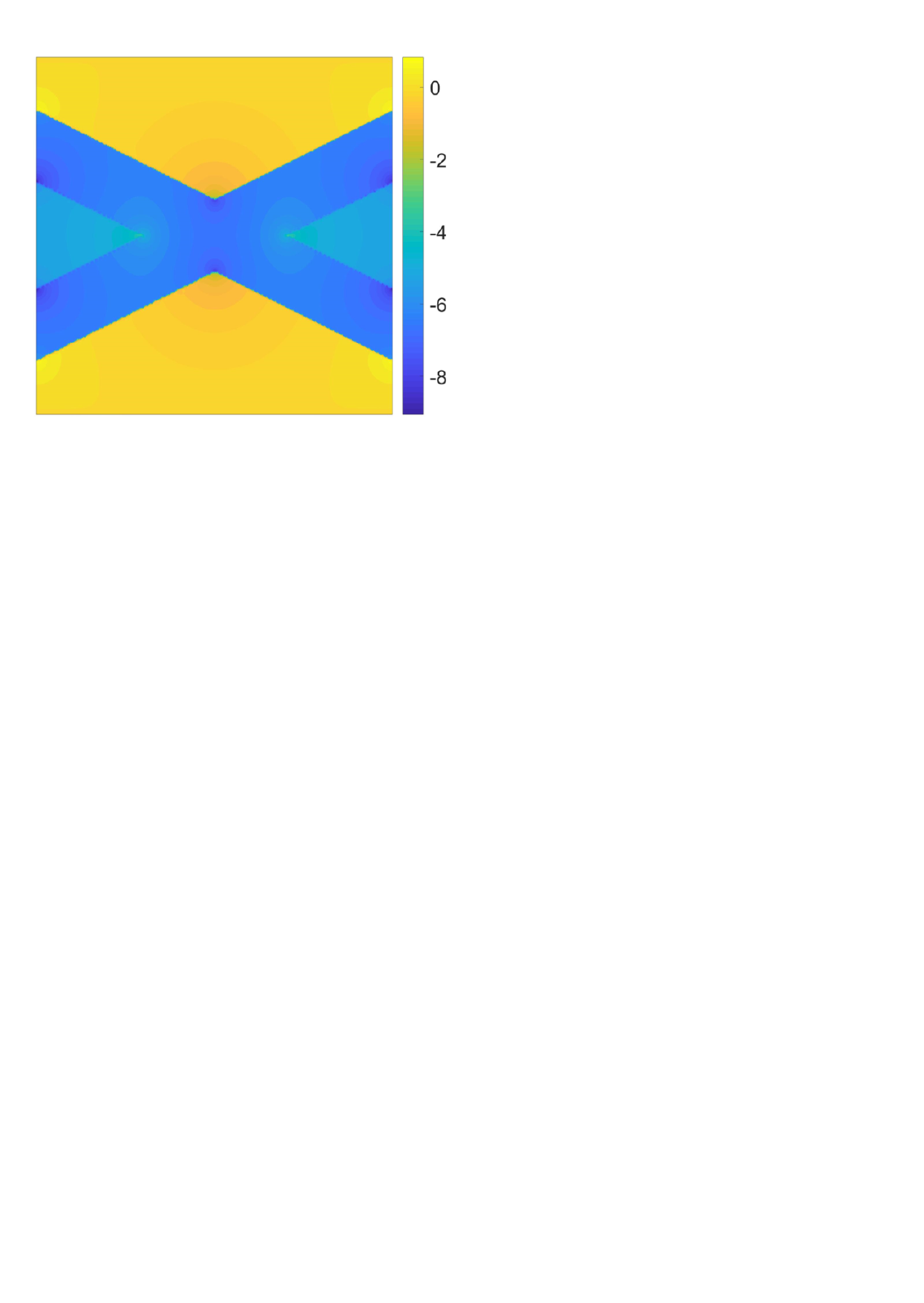}
				\end{center}
			\end{minipage}\vspace*{-6.8cm}
		\end{center}
		\caption{Two intersecting low-permeability streaks (left), numerical pressure (center) and  logarithm of  the  norm  of  the  numerical  velocity  (right).}\label{fig:two:streaks}
	\end{figure}

\begin{figure}[t]
		\begin{center}
			\hspace*{-0.7cm}\begin{minipage}[t]{0.3\textwidth}
				\begin{center}
					\includegraphics[scale=0.32]{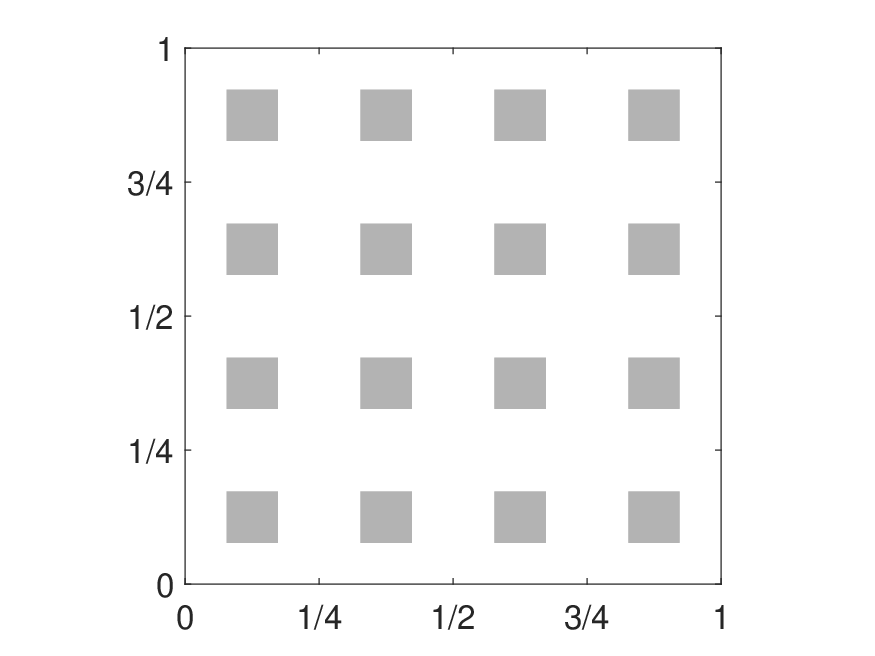}
				\end{center}
			\end{minipage}
			\hspace*{0.3cm}
			\begin{minipage}[t]{0.3\textwidth}
				\begin{center}
					\includegraphics[scale=0.32]{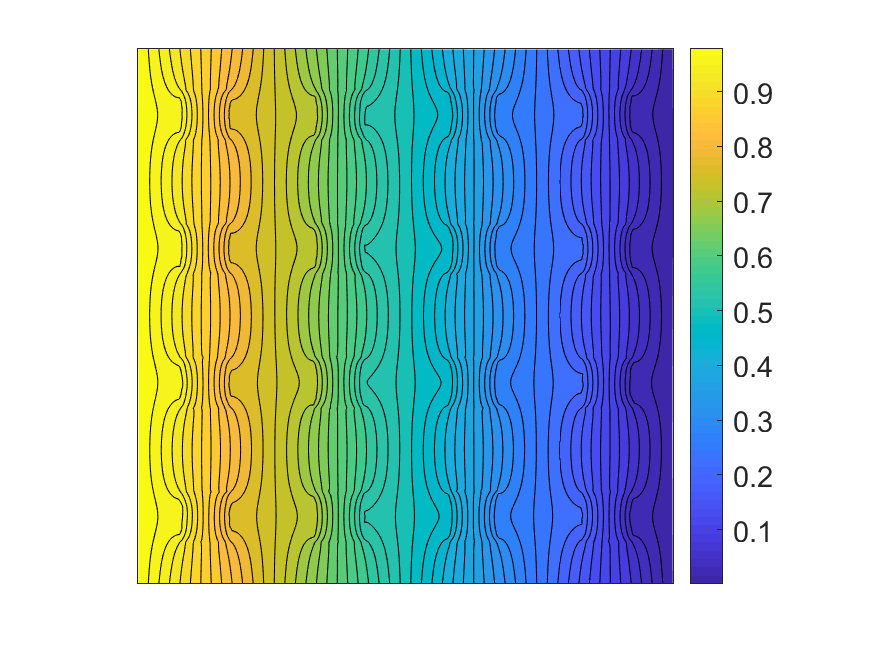}
				\end{center}
			\end{minipage}
			\hspace*{0.7cm}
			\begin{minipage}[t]{0.3\textwidth}
				\begin{center}\vspace*{-0.187cm}
					\includegraphics[scale=0.355]{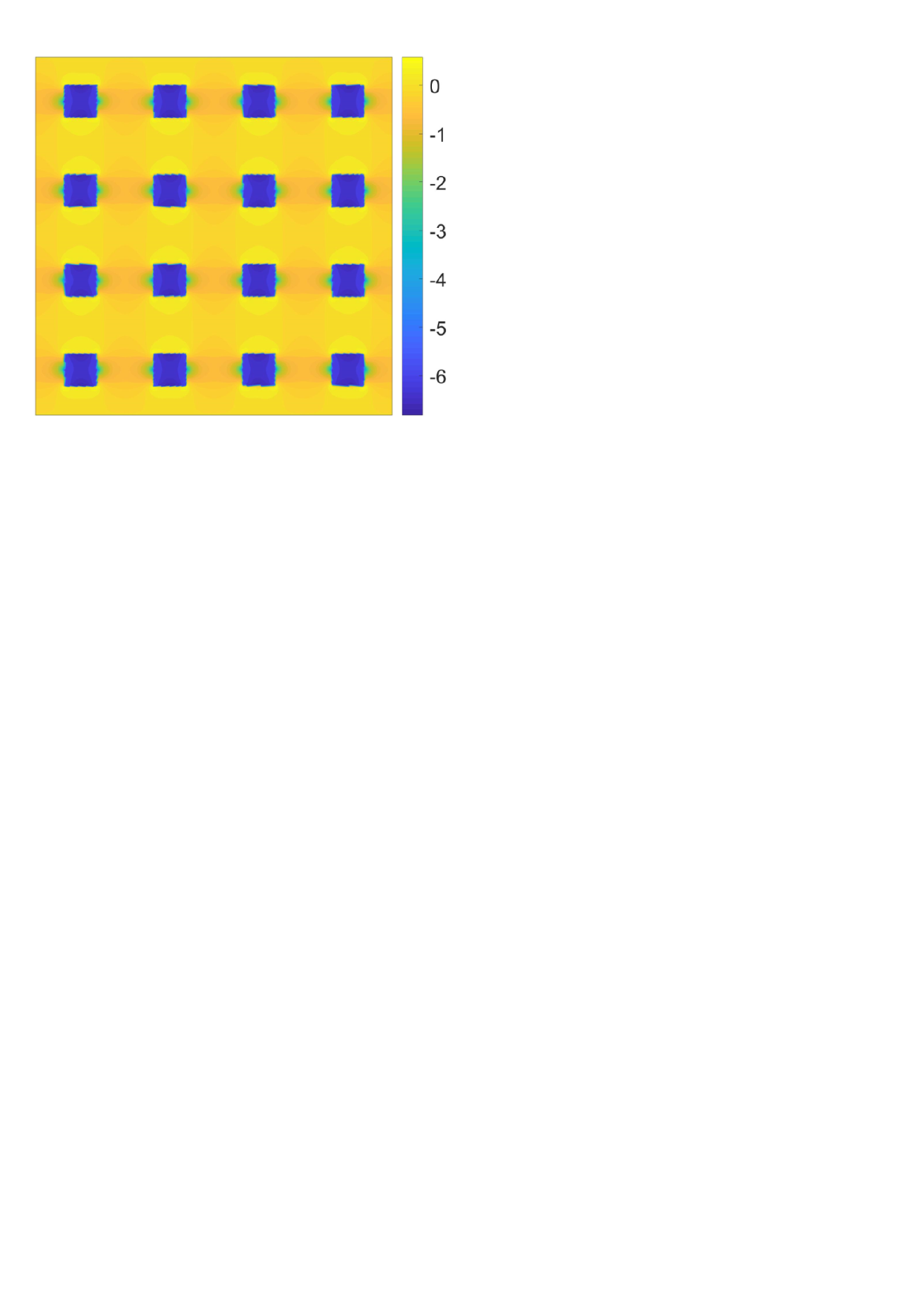}
				\end{center}
			\end{minipage}\vspace*{-6.8cm}
		\end{center}
		\caption{Square-shaped low-permeability inclusions (left), numerical pressure (center) and  logarithm of  the  norm  of  the  numerical  velocity  (right).}\label{fig:square:inclusions}
	\end{figure}

\begin{figure}[t]
		\begin{center}
			\hspace*{-0.7cm}\begin{minipage}[t]{0.3\textwidth}
				\begin{center}
					\includegraphics[scale=0.32]{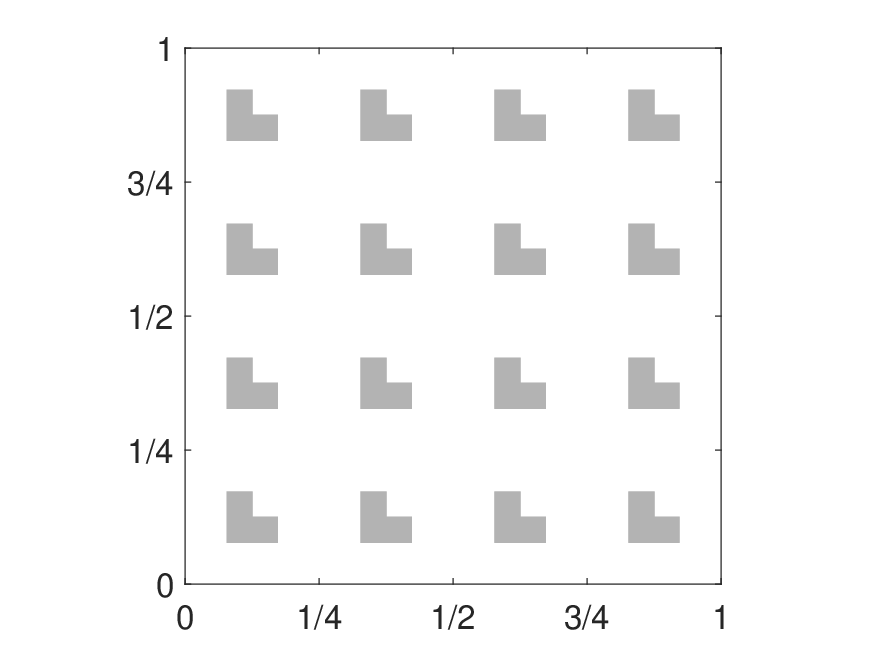}
				\end{center}
			\end{minipage}
			\hspace*{0.3cm}
			\begin{minipage}[t]{0.3\textwidth}
				\begin{center}
					\includegraphics[scale=0.32]{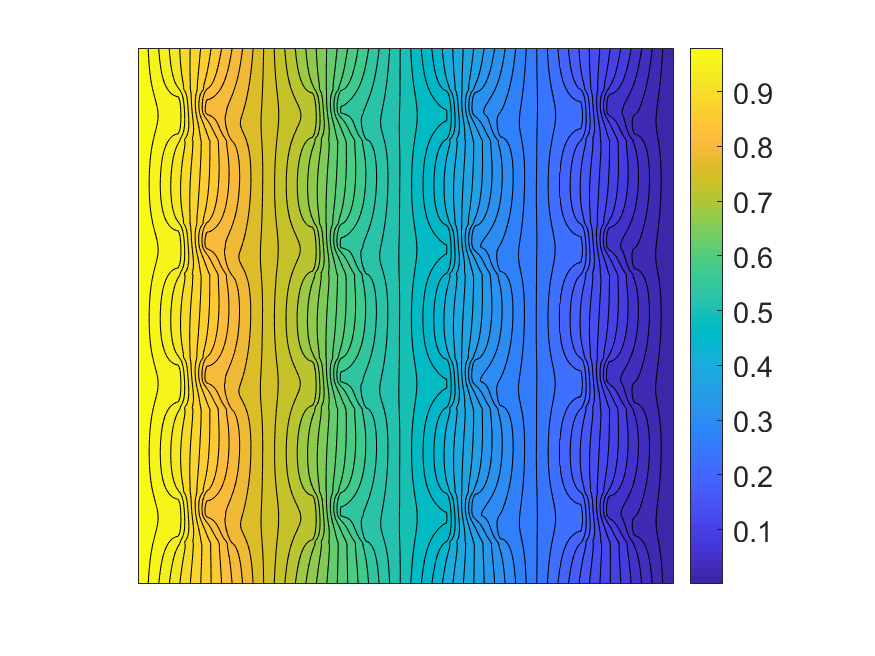}
				\end{center}
			\end{minipage}
				\hspace*{0.7cm}
			\begin{minipage}[t]{0.3\textwidth}
				\begin{center}\vspace*{-0.187cm}
					\includegraphics[scale=0.355]{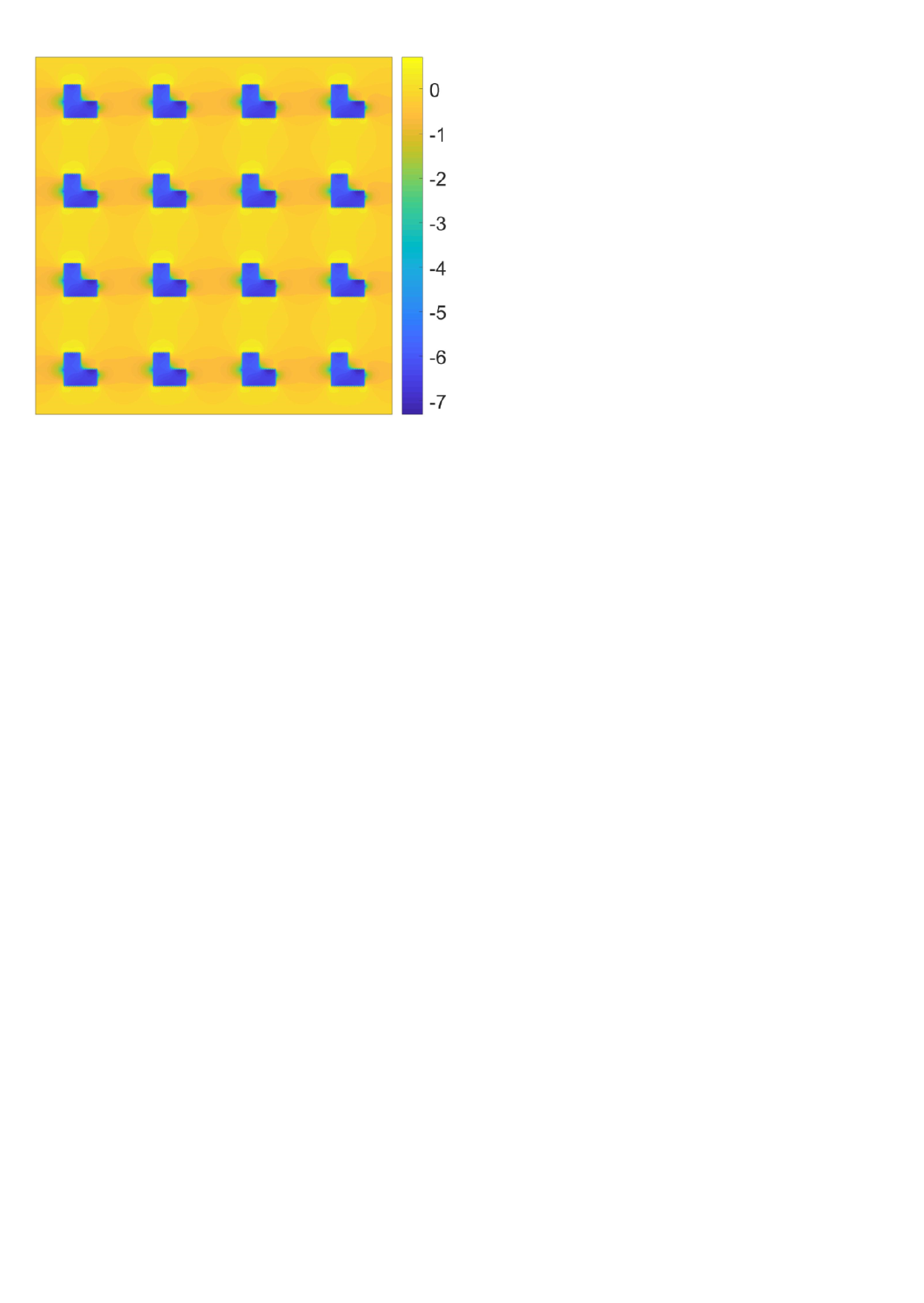}
				\end{center}
			\end{minipage}\vspace*{-6.8cm}
		\end{center}
		\caption{$L$-shaped low-permeability inclusions (left), numerical pressure (center) and  logarithm of  the  norm  of  the  numerical  velocity  (right).}\label{fig:L:inclusions}
	\end{figure}

Next, let us test the behaviour of the multigrid method when considering different configurations of low-permeability regions, and various families of logically rectangular grids which are not aligned with the jumps in the permeability coefficient. 
In particular, we consider two low-permeability intersecting bands, depicted in grey in Figure \ref{fig:two:streaks} (left), and also square-shaped and $L$-shaped low-permeability inclusions, depicted in grey in Figure \ref{fig:square:inclusions} (left) and Figure \ref{fig:L:inclusions} (left), respectively. This type of periodic inclusions were previously considered in \cite{kum:rod:gas:oos:19}. In all the cases, we fix $K=10^{-3} I$ inside the streaks or inclusions, and $K=I$ in the rest of the domain. 
Regarding the spatial meshes, we consider the four families of logically rectangular grids introduced in Section \ref{subsec:3_1} and shown in Figure \ref{fig:meshes}.
	
The central plot in Figures \ref{fig:two:streaks}, \ref{fig:square:inclusions} and \ref{fig:L:inclusions} represents the numerical pressure obtained for the corresponding permeability configuration when considering the \linebreak smooth mesh depicted in Figure \ref{fig:meshes}(a) for $N=320$. Accordingly, the right plot in Figures \ref{fig:two:streaks}, \ref{fig:square:inclusions} and \ref{fig:L:inclusions} represents the logarithm of the norm of the numerical velocity field. In all the cases, the numerical solutions show the expected physical behaviour.
	
Table \ref{table:iter:two:streaks} shows the number of multigrid iterations required to reduce the initial residual by a factor of $10^{10}$ for different values of $h$, different permeability configurations and different families of spatial meshes. In particular, for each permeability distribution, the columns denoted by (a), (b), (c) and (d) correspond to the families of smooth, Kershaw, $h$-perturbed and randomly $h$-perturbed meshes, respectively. Once again, in the case of considering random grids, we show the integer approximation to the average number of iterations required in 100 realizations of the corresponding experiment.
	
	\begin{table}[t]
		\caption{{Number of iterations of the multigrid method when considering different permeability distributions and various families of logically rectangular meshes ($h_0=20^{-1}$).}}
			\label{table:iter:two:streaks}
			\renewcommand{\arraystretch}{1.2}
			\begin{center}
				\begin{tabular}{c|cccc|cccc|cccc}
					\hline\\[-3ex]
					& \multicolumn{4}{|c|}{Two intersecting streaks} & \multicolumn{4}{|c|}{Square-shaped inclusions} & \multicolumn{4}{|c}{$L$-shaped inclusions}\\\hline
					h& \hspace*{0.05cm}(a)\hspace*{0.05cm} & \hspace*{0.05cm}(b)\hspace*{0.05cm} & \hspace*{0.05cm}(c)\hspace*{0.05cm} & (d) & \hspace*{0.05cm}(a)\hspace*{0.05cm} & \hspace*{0.05cm}(b)\hspace*{0.05cm} & \hspace*{0.05cm}(c)\hspace*{0.05cm} & (d) & (a) & (b) & (c) & (d)\\\hline
					$h_0$& 7 & 10 & 7 & 7  & 6 & 9 &7 & 7 & 6 & 9 & 7& 7
					\\
					$h_0/2$& 7  & 9 &7 & 7 & 6 & 9 &6 & 7 & 6 &9  &6 & 7
					\\
					$h_0/2^2$& 7 & 9 &6 & 7 & 6 & 9 & 6& 7 & 7 & 10 &6 & 7
					\\
					$h_0/2^3$&7  & 9 &6 & 6 & 6 & 10 &6 & 7 &7 &10  &7 & 7
					\\
					$h_0/2^4$&7  & 9 &6 & 6 & 6 &10  &6 &7 &8 &11  &7 & 8
					\\\hline
					
				\end{tabular}
			\end{center}
	\end{table}
	
It is worth noting the robustness of the multigrid method for different permeability configurations, even when considering spatial meshes which are not aligned with the jumps on the permeability coefficients.
	

\subsection{A test with a random permeability tensor}\label{sec:test_random}

\begin{figure}[t]
\begin{center}
\begin{tabular}{cc}
\includegraphics[width = 0.37\textwidth]{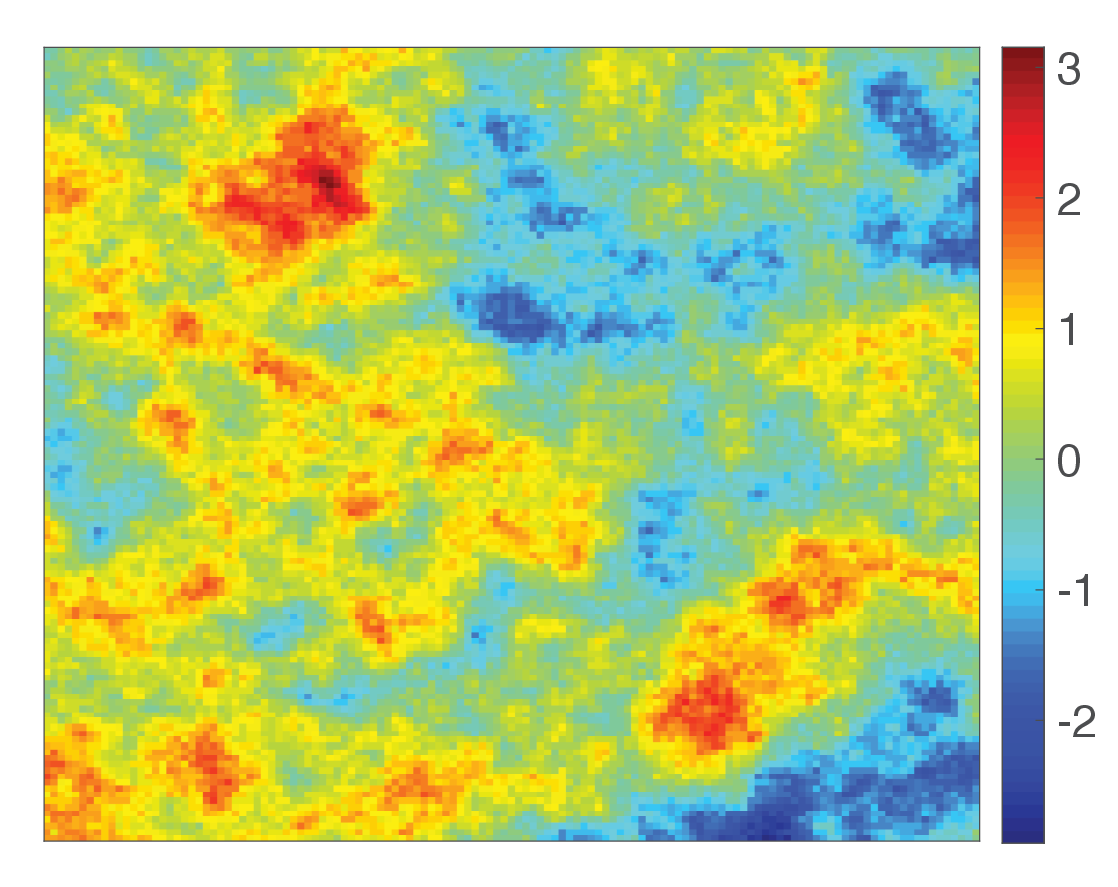}
&
\includegraphics[width = 0.37\textwidth]{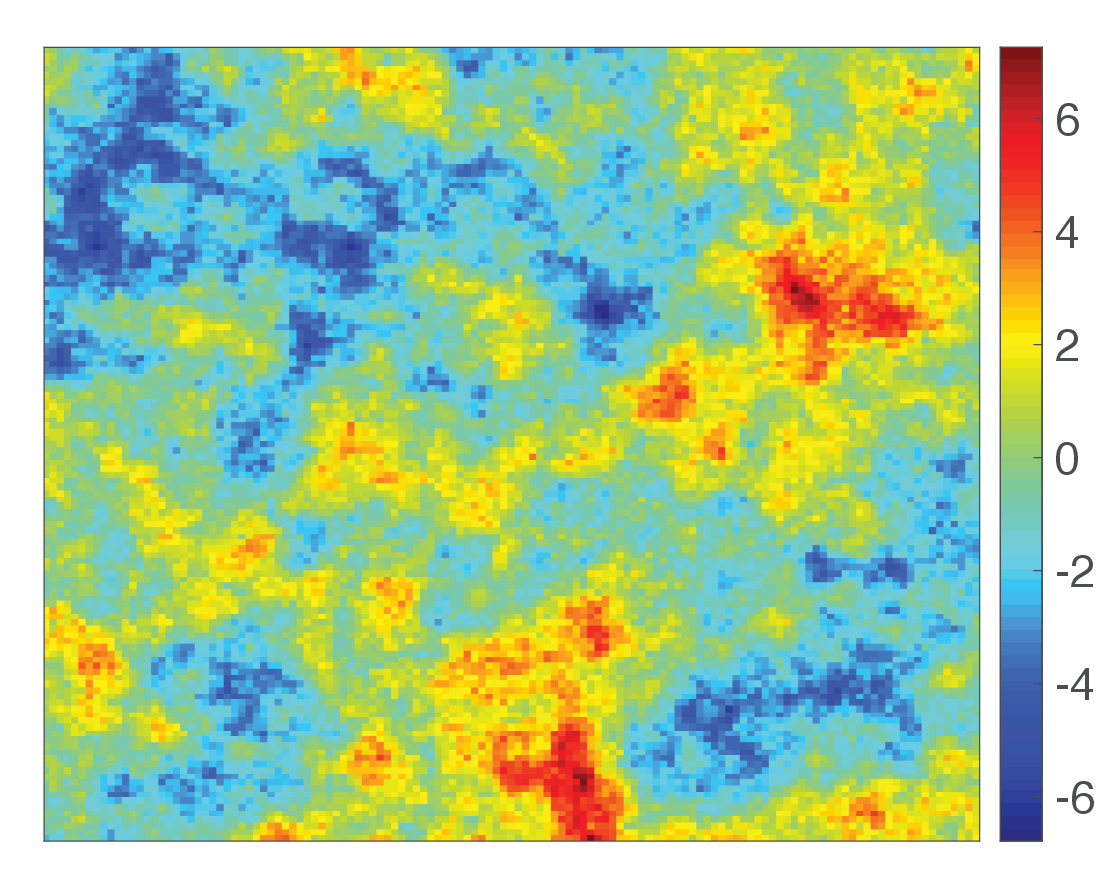}\\
(a) $\Phi_1 = (0.5, 0.3, 1)$ & (b) $\Phi_2 = (0.5, 0.1, 3)$
\end{tabular}
\caption{Logarithm of the permeability field, $\log_{10}k$, generated using two different sets of parameters $\Phi = (\nu_c, \lambda_c, \sigma_c^2)$.}
\label{random_K_pic}
\end{center}
\end{figure}

To conclude, we would like to study the robustness of the proposed multigrid method in the case of considering a heterogeneous random medium. For this purpose, the domain is set to be $\Omega=(0,1)^2$, and Dirichlet conditions and a zero right-hand side are considered. We assume a diagonal permeability tensor $K = kI$ with constant $k$, where $I$ is the identity tensor. A lognormal random field may accurately represent the permeability of a heterogeneous porous medium \cite{Freeze1975}, and therefore we assume that the logarithm of the permeability field, $\log_{10}k$, is modeled by a zero-mean Gaussian random field. In this way, to generate samples of the Gaussian random field, we use the so-called Mat\'ern covariance function \cite{Handcock1994}, given by
\begin{equation}\label{Matern_covariance}
C_{\Phi}(r) = \sigma_c^2 \displaystyle \frac{2^{1-\nu_c}}{\Gamma(\nu_c)}\left(2\sqrt{\nu_c}\frac{r}{\lambda_c}\right)^{\nu_c}K_{\nu}\left(2\sqrt{\nu_c}\frac{r}{\lambda_c}\right).
\end{equation}
Such a covariance function is characterized by a set of parameters $\Phi = (\nu_c, \lambda_c, \sigma_c^2)$, where $\nu_c$ defines the field smoothness, $\lambda_c$ is the correlation length and $\sigma_c^2$ represents the variance. In the previous expression, $r$ is the distance between two points, $\Gamma$ is the gamma function, and $K_{\nu}$ is the modified Bessel function of the second kind. By using the Mat\'ern family of covariance functions, random coefficient fields with different degrees of smoothness can be generated. In particular, we consider two Mat\'ern reference sets of parameters with increasing order of complexity, namely: $\Phi_1 = (0.5, 0.3, 1)$ and $\Phi_2 = (0.5, 0.1, 3)$. In Figure \ref{random_K_pic}, we represent a possible sample of $\log_{10}k$ using each of the two sets of parameters, as an example. We can observe that, for the set $\Phi_2$, the fluctuations of the permeability field are much larger than for $\Phi_1$.

For a fixed set of parameters and a mesh size $h$, we generated $100$ realizations of the permeability field and computed the average number of multigrid iterations necessary to reduce the initial residual by a factor of $10^{-9}$. This experiment was done for both regular and randomly $h$-perturbed meshes and Table \ref{table:iterations_ex3} shows the rounding to the closest natural number of the average number of iterations performed in each case. We can observe from the table that the convergence of the solver is independent of the discretization parameter $h$. For both types of grids, we obtain an average of around $5$-$6$ iterations for both sets of parameters $\Phi_1$ and $\Phi_2$. Summarizing, the performance of the proposed multigrid solver in highly heterogeneous random media is shown to be very satisfactory too.
\begin{table}[t]
		\caption{Average number of multigrid iterations when considering two different sets of parameters $\Phi_1$ and $\Phi_2$ and both regular and randomly $h$-perturbed meshes ($h_0 = 2^{-5}$) in the test with random permeability fields.}
		\label{table:iterations_ex3}
		\renewcommand{\arraystretch}{1.2}
		\begin{center}
			\begin{tabular}{c|cc|cc}
				\hline\\[-3ex]
				Mesh&\multicolumn{2}{c|}{Regular mesh}&\multicolumn{2}{c}{Randomly $h$-perturbed}  \\\hline
				$h$& \hspace*{0.2cm}$\Phi_1$\hspace*{0.2cm} & $\Phi_2$ &\hspace*{0.4cm} $\Phi_1$\hspace*{0.4cm} & $\Phi_2$ 
				\\\hline
				$h_0$&5 & 8 & 5 & 8
				\\
				$h_0/2$&5 & 7 & 5 & 7
				\\				
				$h_0/2^2$& 5 & 6 & 5 & 6  
				\\
				$h_0/2^4$& 5 & 5 & 5 & 5
				\\\hline		
			\end{tabular}
		\end{center}
\end{table}

\section{Conclusions}
\label{sec:conclusions}
The aim of this work is the efficient solution of multipoint flux approximations of the Darcy problem. These methods can be applied on irregular meshes and accurately handle anisotropic discontinuous tensors even with large coefficient jumps. Thus, the search for a solver for MPFA is of great interest.

Here, we have proposed a blackbox-type geometric multigrid method for MPFA on logically rectangular meshes. This is the first time that multigrid is applied to MPFA discretizations. The blackbox methodology makes the solver very easy to apply for the user, and the logically rectangular meshes take advantage of the recent trends in computer architectures that achieve their best performance when structured data can be used. The robustness of the proposed multigrid solver has been shown for several problems with different permeability tensors, including random permeability fields, on a variety of logically rectangular grids. Moreover, the proposed method can be easily parallelized by using grid partitioning techniques.

The main limitation of this solver, related to its applicability to two-dimensional problems, can also be overcome. To this respect, the straightforward extension to three-dimensional problems would imply the use of plane relaxation instead of line-smoothers, which would make the method, however, prohibitively more expensive. Thus, we would propose to apply one iteration of the two-dimensional multigrid presented here on each plane instead of solving it exactly. The extension to three-dimensional problems, however, is a subject of future research.

Finally, another extension for the multigrid strategy would be to consider the higher order MFMFE methods described in \cite{amb:kha:lee:yot:19}. Such methods are based on a new family of enhanced Raviart--Thomas spaces, and also reduce to a cell-centered pressure system if the velocity degrees of freedom are chosen to be the points of tensor-product Gauss--Lobatto quadrature rules. In this context, the blackbox multigrid solver could be naturally applied to the resulting pressure system. Alternatively, we could also design a so-called $p$-multigrid strategy\footnote{This technique was first proposed in \cite{ron:pat:87,mad:mun:88} for spectral element discretizations, and has been subsequently extended to discontinuous Galerkin methods \cite{fid:oli:tod:dar:05} or isogeometric analysis \cite{tie:mol:god:vui:20}.} that takes advantage of different approximation orders along the iterations.

\begin{acknowledgements}
Francisco J. Gaspar has received funding from the European Union's Horizon 2020 Research and Innovation Programme under the Marie Sk{\l}odowska--Curie grant agreement No. 705402, POROSOS. The work of Laura Portero is supported by the Spanish project MTM2016-75139-R (AEI/FEDER, UE) and the Young Researchers Programme 2018 from the Public University of Navarre. Andr\'es Arrar\'as acknowledges support from the Spanish project PGC2018-099536-A-I00 (MCIU/AEI/FEDER, UE) and the Young Researchers Programme 2018 from the Public University of Navarre. The work of Carmen Rodrigo is supported by the Spanish project PGC2018-099536-A-I00 (MCIU/AEI/FEDER, UE) and the DGA (Grupo de referencia APEDIF, ref. E24\_17R). Andr\'es Arrar\'as and Laura Portero gratefully acknowledge the hospitality of the Centrum Wiskunde and Informatica (Amsterdam), where this research was partly carried out.
\end{acknowledgements}

\bibliographystyle{spmpsci}      
\bibliography{Arraras_Gaspar_Portero_Rodrigo_MG_MPFA}   


\end{document}